\documentclass[11pt, a4paper]{article}
\input epsf
\usepackage{amsmath,amsfonts,amssymb,amsthm}
\usepackage[T1]{fontenc}
\newtheorem{thm}{Th\'eor\`eme}

\newtheorem{prop}{Proposition}
\newtheorem{corollary}{Corollaire}
\newcommand\inv{\mathop{\rm inv}}

\newenvironment{remark}{\smallskip\noindent{\bf Remarque. }}{\smallskip}
\newtheorem{lem}{Lemme}

\begin{document}
\title{D\'eveloppements limit\'es et r\'eversion des s\'eries}
\author{R. Bacher et B. Lass}
\date{}

\maketitle
{\it R\'esum\'e: Nous \'etudions quelques 
propri\'et\'es nouvelles li\'ees aux d\'eveloppements
limit\'es et \`a la transform\'ee
de Hankel. Nous les d\'emontrons en utilisant
l'approche combinatoire de la r\'eversion 
des s\'eries et des fractions continues.} 

\section{Introduction}

Le but de cet article est de d\'ecrire quelques interactions entre 
combinatoire et alg\`ebre. Plus pr\'ecis\'ement, nous \'etudions 
certains liens entre la r\'eversion des s\'eries (formelles) et
les matrices de Hankel. Les deux sujets sont classiques : 
la plupart des fonctions importantes, par exemple $\hbox{exp},\ \hbox{sin}$ ou 
$\hbox{tan}$, poss\`edent en effet des fonctions inverses 
($\hbox{log},\ \hbox{arcsin},\ 
\hbox{arctan}$ dans notre cas) et un th\'eor\`eme c\'el\`ebre 
de Lagrange relie le d\'eveloppement en s\'erie d'une fonction analytique
au d\'eveloppement en s\'erie de sa fonction inverse (pour la composition).
Du c\^ot\'e de la combinatoire, inverser des s\'eries g\'en\'eratrices
est une technique standard, par exemple pour la r\'esolution de
probl\`emes d'\'enum\'eration,
voir \cite{GJ}. La formule de Lagrange-B\"urmann est \'egalement 
utile dans l'\'etude de certains aspects des formes modulaires,
voir \cite{MOS}. Les matrices de Hankel
apparaissent naturellement lorsqu'on consid\`ere les moments d'une
mesure de probabilit\'e convenable sur ${\mathbb R}$ et sont
\'etroitement reli\'ees aux polyn\^omes orthogonaux et \`a certaines
fractions continues. Un traitement
combinatoire de ces matrices a \'et\'e donn\'e par exemple par 
Flajolet dans \cite{F} ou par Viennot dans \cite{V} et \cite{V1}. 
Les r\'ef\'erences
\cite{H1} et \cite{Stanl} contiennent \'egalement quelques
informations historiques.

Notre article est organis\'e comme suit. 

Pour la commodit\'e du lecteur, nous rappelons
le th\'eor\`eme de Lagrange
(concernant la r\'eversion des s\'eries) et 
une preuve classique au d\'ebut du chapitre \ref{lagrange}.
Dans le reste du chapitre, 
nous \'enon\c cons notre r\'esultat principal,
un lien entre le th\'eor\`eme de Lagrange et une suite
de d\'eveloppements limit\'es.

Le chapitre \ref{inverse} relie
la suite associ\'ee aux d\'eveloppements limit\'es 
\`a la \lq\lq transform\'ee inverse''.

Le chapitre \ref{chapexples} contient deux exemples
illustrant les r\'esultats \'enonc\'es. 

Le chapitre \ref{interpol}, ind\'ependant du reste, 
d\'ecrit une d\'eformation continue naturelle qui permet d'interpoler
entre l'inversion $\frac{1}{f}$ et la r\'eversion
$(xf)^{\langle-1\rangle}$ d'une s\'erie formelle 
$f=1+x\mathbb C[[x]]$.

Nous discutons ensuite
quelques jolies propri\'et\'es de la matrice de Hankel
associ\'ee \`a la suite obtenue par des d\'eveloppements
limit\'es, au chapitre \ref{hankel}. 

Le chapitre \ref{Lukas} rappelle
une interpr\'etation combinatoire classique qui fait le
lien entre les matrices de Hankel et divers objets combinatoires
(chemins, mots de \L ukasiewicz). Ces ingr\'edients sont
ensuite utilis\'es pour prouver une partie de nos r\'esultats. 
Ce chapitre contient \'egalement 
des preuves succinctes de r\'esultats classiques (\`a
l'exception de la Proposition \ref{indeps1} qui 
est peut-\^etre moins connue) ainsi qu'une digression d\'ecrivant une
action du groupe di\'edral infini sur les mots de \L ukasiewicz.

Le chapitre final contient des r\'esultats concernant les 
d\'eterminants de matrices de Hankel ainsi que les preuves
des r\'esultats non d\'emontr\'es ant\'erieurement.

Dans la suite, nous travaillerons toujours sur un corps
de caract\'eristique z\'ero.

\section{Le th\'eor\`eme d'inversion de Lagrange}\label{lagrange}

L'ensemble des s\'eries formelles du type $p(x) = \sum_{n=1}^{\infty}  
\alpha_nx^n$
telles que $\alpha_1 \ne 0$ (pour un corps de base fix\'e une fois  
pour toutes)
constitue un groupe pour la {\it composition}.
L'inverse $q(x)$ d'une telle s\'erie $p(x)$ est uniquement d\'efini
par l'\'equation $q \circ p (x) = x$,
et d'ailleurs aussi par $p \circ q (x) = x$.
Le passage de $p$ \`a $q$,
qui est ce que nous appelons ici la {\it r\'eversion des s\'eries},
est l'objet d'un th\'eor\`eme c\'el\`ebre de Lagrange,
qui semble avoir son origine historique dans l'article \cite{Lagr}.

Pour la commodit\'e du lecteur, nous indiquons d'abord au th\'eor\`eme
\ref{lagrhist} une d\'emonstration du th\'eor\`eme de Lagrange
sans doute assez proche de l'original,
en suivant le joli article de Henrici \cite{H1}.
Par ailleurs, le chapitre \ref{Lukas} contient une autre
d\'emonstration classique,
bas\'ee sur la combinatoire des mots de \L ukasiewicz.
Il existe de nombreuses autres pr\'esentations, dont
\cite{Brom}, pages 158-161, 
\cite{GJ}, pages 15-18, 
\cite{Goursat}, pages 129-133,
\cite{Polya}, pages 145-149,
\cite{Stanl}, pages 38-39, 
\cite{White}, pages 128-136.

L'ensemble des s\'eries formelles du type $s(x) = \sum_{n=0}^{\infty}  
s_nx^n$
telles que $s_0 \ne 0$
constitue un autre groupe pour la {\it multiplication}.
Le passage de $s(x)$ \`a $1/s(x)$ est ce que nous appelons ici {\it  
l'inversion des s\'eries},
et nous y revenons aux chapitres \ref{inverse} et  \ref{interpol}.
Il faut n\'eanmoins prendre garde au fait que de nombreux auteurs
utilisent le terme \lq\lq inversion\rq\rq \ dans le contexte du 
th\'eor\`eme de Lagrange.

Soit $g(x) = \sum_{n \ge N} \gamma_n x^n$ une s\'erie de Laurent  
formelle \`a une
ind\'etermin\'ee, o\`u $N \in \mathbb Z$ et o\`u les coefficients $ 
\gamma_n$ sont dans
le corps de base~; on pose $\gamma_n = 0$ pour $n < N$.
Pour tout $n \in \mathbb Z$,
nous \'ecrivons $[x^n](g(x))$ le $n$--i\`eme coefficient $\gamma_n$  
de $g(x)$.
Le  th\'eor\`eme de Lagrange, ou de Lagrange--B\"urmann,
\'etablit une relation entre les coefficients $[x^n](q(x))$
et les coefficients
$[x^{n-1}] \left( \big(\frac{x}{p(x)}\big)^n\right)$ pour $q(x)$
la r\'eversion d'une s\'erie 
$p(x)=\sum_{n=1}^\infty\alpha_nx^n$ avec $\alpha_1\not=0$.

\begin{thm} \label{lagrhist} Soient $p(x) = \sum_{n=1}^{\infty}  
\alpha_nx^n$
une s\'erie formelle sans terme constant telle que $\alpha_1 \not= 0$,
et $q(x)$ la s\'erie du m\^eme type telle que
$q \circ p (x) = p \circ q (x) = x$. Alors
$$
n[x^n]\left( q^k(x) \right) \, = \,
k[x^{n-k}]  \left( \Big(\frac{x}{p(x)}\Big)^n\right)
$$
pour tous $n,k \in \mathbb Z$.
\end{thm}

\noindent{\bf D\'emonstration}, d'apr\`es \cite{H1}.
Consid\'erons une s\'erie de Laurent formelle $g(x)$ et la s\'erie de  
Laurent
$$
g \circ q (x) \, = \, \sum_{j \geq N} \gamma_j x^j .
$$
Nous obtenons d'abord
$$
g \circ q \circ p (x) \, = \, g(x) \, = \,
\sum_{j \geq N} \gamma_j p^j(x)
$$
en composant  avec $p(x)$ \`a la source,
puis, pour $n \in \mathbb Z$ arbitraire,
$$
g(x) \frac{p'(x)}{p^{n+1}(x)} \, = \,
\sum_{j \geq N} \gamma_j p^{j-n-1}(x) p'(x)
\leqno{(*)}
$$
en multipliant par $p'(x)/p^{n+1}(x)$.

Nous allons appliquer deux r\`egles de calcul tr\`es simples
pour le calcul des r\'esidus. La premi\`ere concerne les d\'eriv\'ees~:
$[x^{-1}](h'(x)) = 0$ pour toute s\'erie de Laurent formelle $h(x)$~;
en particulier,
$$
[x^{-1}]\left( p^{j-n-1}(x) p'(x) \right)
\, = \,
[x^{-1}]\left( \frac{d}{dx} \frac{p^{j-n}(x)}{j-n} \right)
\, = \,
0
$$
pour $j \ne n$.
La seconde concerne les d\'eriv\'ees logarithmiques~:
$[x^{-1}]\left( \frac{h'(x)}{h(x)} \right) = 1$
pour toute s\'erie de puissance de la forme
$h(x) = \sum_{j=1}^{\infty}\delta_j x^j$ avec $\delta_1 \ne 0$.

En \'egalant les r\'esidus des deux termes de (*), nous trouvons donc
$$
[x^{-1}]\left(  g(x) \frac{p'(x)}{p^{n+1}(x)}  \right)
\, = \,
[x^{-1}] \left( \gamma_n \frac{p'(x)}{p(x)} \right)
\, = \, \gamma_n .
\leqno{(**)}
$$
Comme le r\'esidu de la s\'erie de Laurent
$\frac{1}{n}(g/p^{-n})'(x) = \frac{g'(x)}{np^n(x)} - \frac{g(x)p'(x)} 
{p^{n+1}(x)}$
est nul, nous avons aussi
$$
n [x^n] \left( g \circ q (x) \right) \, = \,
n \gamma_n \, = \,
[x^{-1}]\left( \frac{g'(x)}{p^n(x)} \right) .
\leqno{(***)}
$$
En particulier, lorsque $g(x) = x^k$, nous avons
$$
[x^{-1}]\left( \frac{kx^{k-1}}{p^n(x)} \right) \, = \,
[x^{n-k}]\left( \frac{kx^n}{p^n(x)} \right)
$$
et le th\'eor\`eme r\'esulte de ce cas de l'\'egalit\'e (***).
\hfill $\square$

\noindent
{\bf Remarques.}
(i) Plus g\'en\'eralement, la formule (***) fournit le $n$--i\`eme  
terme de la
s\'erie de Laurent $g \circ q (x)$ pour tout $n \ne 0$,
et la formule (**) pour $n=0$ s'\'ecrit
$$
[x^0](g \circ q (x)) \, = \,
[x^{-1}]\left(  g(x) \frac{p'(x)}{p(x)}  \right) .
$$

(ii) Si les coefficients sont complexes et si le rayon de convergence
de la s\'erie $p(x)$ est strictement positif,
alors il en est de m\^eme de celui de $q(x)$.

(iii) D'un point de vue num\'erique,
la s\'erie de von Neumann
$({\mathcal J}-H)^{\langle -1\rangle}={\mathcal J}+
H\circ({\mathcal J}+
H\circ({\mathcal J}+\dots ))$ permet 
de calculer efficacement la r\'eciproque 
$({\mathcal J}-H)^{\langle -1\rangle}$
d'une perturbation $H$ d'ordre $>1$ de l'identit\'e ${\mathcal J}$
en un nombre quelconque de variables.
Cette formule est l'analogue compositionelle 
de la r\`egle de Horner: $1+H(1+H(1+\dots))=\sum_{n=0}^\infty
H^n$ qui converge vers $\frac{1}{1-H}$ pour $H$ petit.
\medskip

\noindent
{\bf Exemples.}
Le th\'eor\`eme \ref{lagrhist} ne s'applique 
bien au calcul des coefficients de $q$ que 
s'il est facile de d\'eterminer
les coefficients de $\left(\frac{x}{p(x)}\right)^n$.

(i) Si $p(x) = \frac{x}{1+x}$, alors $\frac{x}{p(x)} = 1+x$
et le th\'eor\`eme \ref{lagrhist} 
implique $[x^n](q(x)) = 1$ pour tout $n \ge 1$,
en accord avec les \'egalit\'es
$q(x) = \frac{x}{1-x} = \sum_{n=1}^\infty x^n$.
(Notons que les deux s\'eries $p(x)$ et $q(x)$
convergent dans le disque unit\'e.)

(ii) Si $p(x) = xe^{-x}$, on obtient sans peine
$q(x) = \sum_{n=1}^{\infty} \frac{n^{n-1}}{n!}x^n$.
(Notons que, dans ce cas, le rayon de convergence de $p(x)$ est infini
et celui de $q$ est $e^{-1}=\hbox{lim}_{n\rightarrow\infty}
\frac{n^{n-1}}{n!}\ \frac{(n+1)!}{(n+1)^{n}}$.)

De mani\`ere analogue, $p(x)=xe^{-x^2}$ donne $q(x)=\sum_{j=0}^\infty
\frac{(2j+1)^{j-1}}{j!}x^{2j+1}$. (Le rayon de convergence de $p(x)$
est de nouveau infini tandis que la s\'erie de $q(x)$ converge
absolument pour $\vert x\vert<1/\sqrt{2 e}$.)
 
(iv) Comme d\'eja mentionn\'e,
la formule d'inversion de Lagrange-B\"urmann n'est que rarement
utile pour la r\'eversion d'une s\'erie formelle. Des m\'ethodes 
diff\'erentes sont g\'en\'eralement beaucoup plus simples. Un tel exemple
est la fonction $p(x^2) = (\sin x)^2$ ; voir la 
page 130 de \cite{White}. Nous ne savons pas utiliser 
la formule de Lagrange pour prouver que la r\'eversion de
$p(x)$ est donn\'ee par la fonction hyperg\'eom\'etrique
$q(x)=\sum_{j=1}^\infty 2^{2j-1}\frac{x^j}{j^2{2j\choose j}}$.
Un calcul facile montre cependant que $q(x)$ est une solution 
(formelle) de l'\'equation diff\'erentielle 
$$(x^2-x)y''+(x-\frac{1}{2})y'+\frac{1}{2}=0\ .$$
En d\'erivant $q(\sin^2(\sqrt{z}))=z$ par rapport \`a $z$
et en posant $x=\sin^2\sqrt{z}$, nous trouvons
$$q'(x)=\frac{\mathrm{arcsin}\sqrt{x}}{\sqrt{x(1-x)}}$$
et ensuite 
$$q''(x)=\frac{1}{2x(1-x)}\left(1-(1-2x)
\frac{\mathrm{arcsin}\sqrt{x}}{\sqrt{x(1-x)}}\right)\ .$$
Ceci montre que $q(x)$ est \'egalement
solution de l'\'equation diff\'erentielle ci-dessus. 
Un d\'eveloppement \`a l'ordre deux des deux s\'eries permet de 
conclure.
\bigskip

Les deux th\'eor\`emes qui suivent fournissent d'autres paires
du type $(p(x),q(x))$. 
L'aspect peut--\^etre original de notre exposition
consiste \`a faire jouer un r\^ole important aux poly\-n\^omes $P_j$
(et plus bas aux polyn\^omes $Q_j$),
que nous voyons comme des {\it d\'eveloppements limit\'es}
des s\'eries correspondantes.
Si $s(x) = \sum_{j=0}^{\infty} s_j x^j$ est une s\'erie enti\`ere
et $k$ un entier positif, nous notons
$\lfloor s(x) \rfloor_k = s_0 + s_1x + \cdots + s_{k-1}x^{k-1}$
son d\'eveloppement limit\'e \`a l'ordre $k-1$.

Consid\'erons une s\'erie formelle $s(x) = \sum_{j=0}^{\infty}s_jx^j$
telle que $s_0 \ne 0$. D\'efinissons successivement
\begin{itemize}
\item
les polyn\^omes
$$
P_1(x) =  s_0,\quad
P_2(x) =  s_0^2 + s_0s_1 x,\quad
P_3(x) =  s_0^3 + (s_0^2s_1+s_0s_1s_0)x + (s_0^2s_2+s_0s_1^2)x^2,\,\dots
$$
d\'efinis r\'ecursivement par $P_k(x) \, = \, \lfloor P_{k-1}(x)s(x)  
\rfloor_k$,
\item
les constantes $Q_n(0)= [x^{n-1}]P_n(x)$, $n \ge 1$, obtenues 
en consid\'erant les coefficients de plus haut degr\'es 
dans les polyn\^omes $P_1(x),P_2(x),\,\dots$, o\`u $P_n(x)$ est
consid\'er\'e comme \'etant de degr\'e $n-1$,
\item
la s\'erie g\'en\'eratrice
$$
q(t) \, = \, \sum_{n=1}^{\infty} Q_n(0)t^n
$$
des nombres $Q_n(0)$.
\end{itemize}

\begin{thm} \label{invlagrange}
La s\'erie formelle $q(t)$
associ\'ee comme ci--dessus \`a $s(x) = \sum_{j=0}^{\infty}s_jx^j$
v\'erifie
$$
q(t) \, = \, t s(q(t)) .
$$
\end{thm}

En posant $p(x) = x/s(x)$, on retrouve deux s\'eries $p(x), q(x)$
telles que
$$
p(q(t)) \, = \, \frac{q(t)}{s(q(t))} \, = \,
\frac{ ts(q(t))}{s(q(t))} \, = \, t .
$$
Avec ces nouvelles notations, le th\'eor\`eme \ref{lagrhist} 
s'\'ecrit comme suit.

\begin{thm} \label{lagrangerev}
Si  $q(t) = ts(q(t))$, alors 
$$
q(t)^{k+1} =  \sum_{n=k+1}^\infty t^n  \frac{k+1}{n} [x^{n-k-1}] 
\left(s(x)^n\right)$$
pour tout $k \in \{0,1,2,\,\dots\}$ et, en particulier
$$
q(t) =   \sum_{n=1}^\infty  \frac{t^n }{n} [x^{n-1}] \left(s(x)^n\right).
$$
\end{thm}

Nous offrons au chapitre \ref{Lukas} une autre preuve des 
th\'eor\`emes \ref{invlagrange} et \ref{lagrangerev}
(bien que ce dernier ne soit rien d'autre qu'une reformulation
du th\'eor\`eme \ref{lagrhist}). Cette preuve, de nature
combinatoire, n'est pas nouvelle. Elle
consiste \`a interpr\'eter les mots de 
\L ukasiewicz comme des arbres plans enracin\'es.


\section{La transform\'ee inverse}\label{inverse}

Le but de ce chapitre est de d\'ecrire quelques aspects du groupe
multiplicatif constitu\'e des s\'eries formelles du type
$\sum_{n=0}^\infty s_nx^n$ avec $s_0\not= 0$. Rappelons qu'une
telle s\'erie d\'efinit une suite de polyn\^omes 
$P_1(x)=s_0,\,\dots,\,P_k(x)=\lfloor P_{k-1}(x)s(x)\rfloor_k,\,\dots,$
o\`u $P_k(x)$ est le d\'eveloppement limit\'e \`a l'ordre $k-1$
de la s\'erie formelle $P_{k-1}(x)s(x)$. Introduisons maintenant les
polyn\^omes  miroir
$Q_n(x)=x^{n-1}P_n(1/x)$ et d\'esignons par $q(t)=\sum_{n=1}^\infty Q_n(0)t^n$
la s\'erie g\'en\'eratrice associ\'ee \`a la suite des 
\'evaluations $Q_1(0),\,Q_2(0),\,\dots$.

Le r\'esultat suivant exprime
la s\'erie g\'en\'eratrice compl\`ete $Q(x)=
\sum_{n=1}^\infty Q_n(x)t^n$ en fonction de $q(t)$:

\begin{thm} \label{thmA} On a 
$$\sum_{n=1}^\infty Q_n(x)t^n=\frac{q(t)}{1-xq(t)}\ .$$
\end{thm}

Ce th\'eor\`eme sera d\'emontr\'e au chapitre \ref{Lukas}.
La preuve consiste \`a
identifier les mon\^omes contribuant
aux coefficients de $q(t)$ avec les mots de \L ukasiewicz.

Nous d\'ecrivons maintenant une interpr\'etation en termes 
de \lq\lq transform\'ee inverse continue'' de cette \'egalit\'e.
Cette interpr\'etation sugg\`ere une jolie 
propri\'et\'e des transform\'ees de Hankel (d\'ecal\'ees) de la s\'erie
$Q_1(x),Q_2(x),\,\dots$  qui sera \'enonc\'ee au chapitre
\ref{hankel} et qui constitue le r\'esultat principal dans
cet article.

Soit $a(t) = a_0 + a_1 t + a_2 t^2 + a_3 t^3 +a_4 t^4 +\cdots$ une
s\'erie g\'en\'eratrice. Introduisons l'application
$I[a(t)] = a(t)/\bigl(1+ta(t)\bigr)$ appell\'ee la 
\emph{transform\'ee inverse}
puisque $\bigl(1+ta(t)\bigr)\bigl(1-tI[a(t)]\bigr) = 1$. Par it\'eration, on obtient
$I^x[a(t)] =  a(t)/\bigl(1+xta(t)\bigr)$ ce qui 
permet d'interpoler les it\'er\'ees 
$$I^x(a(t))=\sum_{k=0}^\infty I_k(x)\ t^k$$
de la transform{\'e}e inverse. Le $k-$i\`eme terme $I_k(x)$ de la suite
$$I_0(x)=a_0,\ I_1(x)=a_1-a_0^2x,\ I_2(x)=a_2-2a_0a_1x+a_0^3x^2
,\,\dots$$
est alors un polyn\^ome de degr\'e $k$ en $x$. 

Posons $a(t)=\frac{q(t)}{t}=\frac{1}{t}\sum_{n=1}^\infty Q_n(0)t^n$. Le
th\'eor\`eme \ref{thmA} s'\'enonce aussi sous la forme
$$I^x\left[\frac{q(t)}{t}\right]=\frac{1}{t}\frac{q(t)}{1+xq(t)}=\frac{1}{t}
\sum_{n=1}^\infty Q_n(-x)t^n\ .$$
Autrement dit, on a $t\ I^x[t^{-1}Q(0)]=Q(-x)$ pour
$Q(x)=\sum_{n=1}^\infty Q_n(x)t^n$.

\begin{remark} On aurait tout aussi bien pu d\'efinir la transform\'ee
de Hankel de $a(t)$ par la formule $\tilde I[a(t)]=a(t)/(1-ta(t))
=-I[-a(t)]$.
\end{remark}

\begin{remark}
Un ph{\'e}nom{\`e}ne similaire d'interpolation continue se produit
{\'e}galement pour la composition it{\'e}r{\'e}e $f^{\circ k}=f\circ f\circ
\dots\circ f$ d'une s{\'e}rie formelle $f(t)=t+\sum_{i=2}^\infty a_i t^i$
dont le d{\'e}veloppement {\`a} l'ordre $1$ est l'identit{\'e}; ceci se 
g{\'e}n{\'e}ralise d'ailleurs facilement {\`a} un d$-$uplet de s{\'e}ries formelles 
$F(t_1,\,\dots,t_d)=(f_1(t_1,\,\dots,t_d),\dots,f_d(t_1,\,\dots,t_d))$. 
Il existe alors une suite
$$\begin{array}{l}
C_1(x)=1,\ C_2(x)=a_2x,\ C_3(x)=(a_2^2(x-1)+a_3)x,\\
\quad C_4(x)=(((2x-3)a_2^3+5a_2a_3)(x-1)+2a_4)x/2,\,\dots\end{array}$$
avec $C_n(x)$ un polyn\^ome de degr{\'e} $\leq n-1$ en $x$ tel que
$f^{\circ x}(t)=\sum_{i=1}^\infty C_i(x)t^i$.

Pour le d\'emontrer on peut consid{\'e}rer la diff{\'e}rence finie
$$C_n(k+1)-C_n(k)=\hbox{coefficient de }t^n\hbox{ dans }
\sum_{i=2}^{\infty}a_i\left(\sum_{j=1}^{\infty}C_j(k)t^j
\right)^i$$
qui est un polyn\^ome de degr{\'e} au plus $n-2$ en 
$k$ (par r{\'e}currence sur 
$n$). On peut {\'e}galement le d{\'e}duire en utilisant un isomorphisme
de mono\"{\i}de entre le mono\"{\i}de des s\'eries 
formelles sans terme constant
(avec la composition des s\'eries comme produit)
et un groupe de matrices triangulaires sup\'erieures.
Un tel isomorphisme peut \^etre donn\'e par
$$\sum_{n=1}^\infty a_n x^n\longmapsto
\left(\begin{array}{ccccc}
a_{1,1}&a_{1,2}&a_{1,3}&a_{1,4}&\dots\\
0&a_{2,2}&a_{2,3}&a_{2,4}\\
0&0&a_{3,3}&a_{3,4}&\dots\\
 & &       &\ddots\end{array}\right)$$
o\`u $\sum_{j=k}^\infty a_{k,j}x^j=\left(\sum_{n=1}^\infty a_n x^n
\right)^k$ ; voir par exemple le
th\'eor\`eme 1.7a dans \cite{H}.
\end{remark}


\section{Exemples}\label{chapexples}

Revenons aux exemples (i) et (ii) du chapitre \ref{lagrange}.

\noindent
{\bf Exemple trivial.} Consid\'erons la s\'erie formelle $s$ d\'efinie par
le polyn\^ome $1+x$. On v\'erifie facilement que
$P_n(x)=Q_n(x)=(1+x)^{n-1}$ et
$q(t)=\sum_{n=1}^\infty t^n=\frac{t}{1-t}$.
Les th\'eor\`emes 1, 2 et 3 se r\'eduisent alors \`a des identit\'es 
triviales et au th\'eor\`eme binomial, \`a savoir
$$\frac{t}{1-t(1+x)}=\frac{t/(1-t)}{1-xt/(1-t)},\ 
\frac{t/(1-t)}{1+t/(1-t)}=t$$
et
$$\left(\frac{t}{1-t}\right)^{k+1}=
\sum_{n=k+1}^\infty \frac{k+1}{n}{n\choose k+1}t^n=\sum_{n=k}^\infty
{n\choose k}t^{n+1}\ .$$

La transform\'ee inverse de $q(t)/t$ est donn\'ee par
$I\left(\frac{q(t)}{t}\right)=\frac{1/(1-t)}{1+t/(1-t)}=1$ et nous avons
$$I^x\left(\frac{q(t)}{t}\right)=
\frac{1/(1-t)}{1+xt/(1-t)}=\frac{1}{1-t(1-x)}
=\sum_{n=1}^\infty
\left(t(1-x)\right)^{n-1}=\frac{1}{t}\sum_{n=1}^\infty
Q_n(-x)t^n$$
en accord avec les r\'esultats du chapitre \ref{inverse}.

\noindent
{\bf L'exemple de l'exponentielle.}
Pour la s\'erie $s(x)=e^x=\sum_{n=0}^\infty \frac{x^n}{n!}$
d\'efinissant l'exponentielle, nous avons
$$P_n(x)=\frac{1}{n}
\sum_{j=0}^n(n-j)\frac{(nx)^j}{j!}.$$
En effet, cette formule donne bien $P_1=1$ et le calcul
\begin{eqnarray*}
&\frac{1}{n}\sum_{j=0}^k(n-j)\frac{n^j}{j!\ (k-j)!}\\
=&\frac{1}{n\ k!}\big(n(n+1)^k-nk(n+1)^{k-1}\big)\\
=&\frac{1}{n+1}(n+1-k)\frac{(n+1)^k}{k!}
\end{eqnarray*}
du coefficient $x^k,0\leq k\leq n$ dans $P_n(x)e^x$ la montre
par r\'ecurrence. Nous obtenons ainsi
$$Q_n(x)=\sum_{j=1}^n j\ n^{n-1-j}\frac{x^{j-1}}{(n-j)!}$$
et
$$q(t)=\sum_{n=1}^\infty Q_n(0)t^n=\sum_{n=1}^\infty n^{n-2}
\frac{t^n}{(n-1)!}=\sum_{n=1}^\infty \frac{1}{n}
\frac{(nt)^n}{n!},$$
en accord avec le th\'eor\`eme \ref{lagrangerev}. Le th\'eor\`eme
\ref{thmA} implique les \'egalit\'es
$$(k+1)(n-k)n^{n-2-k}=k\sum_{m=k}^{n-1}{n-k\choose n-m}m^{m-1-k}
(n-m)^{n-m-1}$$
pour tous les entiers $n,k$ tels que $n>k>1$.
Pour finir, mentionnons la jolie \'evaluation
$$P_n(1)=Q_n(1)=\frac{1}{n}\sum_{j=0}^n(n-j)\frac{n^j}{j!}=
\sum_{j=0}^n\frac{n^j}{j!}-\sum_{j=0}^{n-1}\frac{n^j}{j!}=
\frac{n^n}{n!}.$$

\section{Interpolation entre inversion et r\'eversion 
d'une s\'erie formelle}\label{interpol}

L'anneau ${\mathbb C}[[x]]$ des s\'eries formelles est un anneau
commutatif local dont l'id\'eal maximal ${\mathfrak m}=x{\mathbb
  C}[[x]]$ est l'ensemble des s\'eries formelles sans terme constant.
Notons $${\mathcal U}={\mathbb C}[[x]]\setminus {\mathfrak m}
={\mathbb C}^*+{\mathfrak m}$$
le groupe multiplicatif form\'e des \'el\'ements inversibles de 
${\mathbb C}[[x]]$ et $\mathcal{SU}=1+{\mathfrak m}\subset 
{\mathcal U}$ le
  sous-groupe des s\'eries formelles de
coefficient constant $1$.
Notons $${\mathcal D}={\mathfrak m}\setminus {\mathfrak m}^2=
\{\sum_{j=1}^\infty \alpha_j x^j\in{\mathbb C}[[x]]\vert
\alpha_1\not=0
\}$$ le groupe non-commutatif des s\'eries formelles
pour la composition. On a ${\mathcal D}=x{\mathcal U}$ en tant
qu'ensemble et $\mathcal{SD}=x\mathcal{SU}=x+{\mathfrak m}^2$ 
peut \^etre interpr\'et\'e comme le sous-groupe des 
``diff\'eomorphismes locaux formels tangents  
\`a l'identit\'e en $0$''. 

Le but de ce chapitre est de d\'ecrire une d\'eformation 
naturelle continue 
(qui est holomorphe pour des s\'eries holomorphes) entre le groupe
multiplicatif commutatif $\mathcal{SU}$ et le 
groupe non-commutatif $\mathcal{SD}$ (identifi\'e \`a $\mathcal{SU}$ 
via la bijection ensembliste $A\longrightarrow xA$ de 
$\mathcal{SU}$ sur $\mathcal{SD}$). 

L'action naturelle $\alpha\cdot A=A\circ \alpha$
de $\alpha\in{\mathcal D}$ sur un \'el\'ement $A\in{\mathbb C}[[x]]$
agit par automorphismes sur ${\mathcal U}$ et 
$\mathcal{SU}$ et on peut donc former
le produit semi-direct $\mathcal{I}={\mathcal U}\rtimes
{\mathcal D}$ qui est un groupe pour la loi de composition
$$(A,\alpha)(B,\beta)=(C,\gamma)=(A(B\circ\alpha),
\beta\circ \alpha)$$
o\`u $C=A(B\circ\alpha)$ est le produit de la s\'erie 
$A$ avec la s\'erie $B\circ \alpha$. 
L'\'el\'ement inverse $(A,\alpha)^{-1}$ de
$(A,\alpha)$ est donn\'e par
$$(A,\alpha)^{-1}=\left(\frac{1}{A\circ\alpha^{\langle
    -1\rangle}},\alpha^{\langle -1\rangle}\right)$$
o\`u la r\'eversion (ou s\'erie r\'eciproque) $\alpha^{\langle
  -1\rangle}$ de $\alpha\in{\mathcal D}$ est d\'efinie par l'identit\'e
$\alpha\circ\alpha^{\langle -1\rangle}=
\alpha^{\langle -1\rangle}\circ\alpha=x$.
On a les homomorphismes $A\longmapsto (A,x)$ et $(A,\alpha)\longmapsto
\alpha$ (avec section $\alpha\longmapsto (1,\alpha)$) 
provenant de la suite exacte scind\'ee \'evidente
$$0\longrightarrow {\mathcal U}\longrightarrow \mathcal I=
{\mathcal U}\rtimes
  {\mathcal D}\longrightarrow {\mathcal D}\longrightarrow 1\ .$$

Notons $\mathcal{SI}=\mathcal{SU}\rtimes \mathcal{SD}$ le noyau 
$\mathrm{Ker}(\psi)$ de l'homomorphisme de groupes 
$\psi:\mathcal I\longrightarrow 
{\mathbb C}^*\times {\mathbb C}^*$ d\'efinie par 
$\psi(A,\alpha)=(A(0),\alpha'(0))$.

\begin{remark} (i) Le groupe $\mathcal I$ peut se g\'en\'eraliser
  facilement en consid\'erant le produit semi-direct $U\rtimes
D$ o\`u $U$ est un groupe de germes de fonctions
inversibles au voisinage d'un point $P\in X$ avec $X$ un espace
topologique et o\`u $D$ est un groupe de germes
d'hom\'eomorphismes avec point fixe $X$. En particulier, on peut, au 
moins formellement,
remplacer le groupe multiplicatif ${\mathcal U}$ par
le groupe multiplicatif des s\'eries de Laurent non-nulles.

\ \ (ii) Le noyau $\mathcal{SI}=\mathrm{Ker}(\psi)=\mathcal{SU}\rtimes 
\mathcal{SD}$ 
est contractile pour une topologie raisonnable sur ${\mathbb C}[[x]]$ 
(obtenu par exemple en consid\'erant la convergence 
coefficient par coefficient). 
On a donc $\pi_1(\mathcal I)=\pi_1({\mathbb C}^*\times {\mathbb C}^*)=
{\mathbb Z}^2$ pour le groupe fondamental $\pi_1(\mathcal I)$ 
et on peut consid\'erer l'extension centrale
$$0\longrightarrow{\mathbb Z}^2\longrightarrow \tilde{\mathcal I}
\longrightarrow\mathcal I\longrightarrow 1$$
d\'efinissant le rev\^etement universel
$\tilde{\mathcal I}$ de $\mathcal I$, obtenu
en relevant l'extr\'emit\'e des chemins continus issus du neutre
$(1,x)\in \mathcal I$ ou, de mani\`ere \'equivalente,
en consid\'erant des rel\`evements r\'eels des arguments de
$A(0),\alpha'(0)\in\mathbb C^*$ pour $(A,\alpha)\in\mathcal I$.

\ \ (iii) Le groupe abstrait $\mathcal I$ est isomorphe \`a un
``sous-groupe de Lie'' dans les matrices triangulaires inf\'erieures
infinies, voir \cite{Baint}.
\end{remark}

Pour la description de l'interpolation entre le groupe multiplicatif
${\mathcal U}$ et le groupe
non-commutatif ${\mathcal D}$ il faut 
soit se restreindre 
au sous-groupe $\mathcal{SI}=\mathcal{SU}\rtimes \mathcal{SD}=
\mathrm{Ker}(\psi)$ 
qu'on pourrait appeller le {\it groupe d'interpolation sp\'ecial}
soit travailler dans un groupe interm\'ediaire entre $\mathcal I$ et
son rev\^etement universel $\tilde{\mathcal I}$. 
Nous allons d\'ecrire en d\'etail le
premier cas. Le deuxi\`eme cas est trait\'e bri\`evement dans
\cite{Baint}.

Pour $\tau\in{\mathbb C}$, introduisons le sous-ensemble
$$\mathcal{SG}(\tau)=\{(A,xA^\tau)\ \vert\ A\in \mathcal{SU}=
1+x{\mathbb  C}
[[x]]\}\subset \mathcal{SI}$$ o\`u l'on choisit l'unique d\'etermination
``continue'' du logarithme des s\'eries
formelles de mani\`ere \`a avoir 
$A^\tau=e^{\tau\log A}\in \mathcal{SU}=
1+x{\mathbb C}[[x]]$ 
pour $A\in \mathcal{SU}$.

\begin{prop} \label{Ilambda}
(i) L'ensemble $\mathcal{SG}(\tau)$ est un sous-groupe pour tout
  $\tau\in{\mathbb C}$.

\ \ (ii) Le groupe $\mathcal{SG}(0)$ est isomorphe au 
groupe commutatif $\mathcal{SU}$.

\ \ (iii) Pour $\tau\not=0$ les groupes $\mathcal{SG}(\tau)$ 
sont tous isomorphes au groupe non-commutatif
$\mathcal{SD}$. Un isomorphisme est donn\'e par 
$\alpha\longmapsto
\left(\left(\frac{\alpha}{x}\right)^{1/\tau},\alpha\right)
\in\mathcal{SG}(\tau)$ pour $\alpha\in\mathcal{SD}$.
\end{prop}

\begin{corollary} Pour $\tau\in [0,1]$,
l'application
$$\tau\longmapsto
\frac{1}{A\circ(xA^\tau)^{\langle -1\rangle}}$$
est une d\'eformation continue reliant l'inverse multiplicatif
$\frac{1}{A}$ de $A\in\mathcal{SU}$ \`a la s\'erie r\'eciproque 
$\frac{x}{A\circ(xA)^{\langle -1\rangle}}=
(xA)^{\langle -1\rangle}$ de $(xA)\in\mathcal{SD}$.
\end{corollary}

\noindent
{\bf Id\'ee de la preuve de la proposition \ref{Ilambda}} 
L'assertion (ii) est \'evidente.

Un petit calcul montre que l'application $(A,\alpha)\longmapsto 
\left(A\left(\frac{\alpha}{x}\right)^{\lambda},\alpha\right)$
est un automorphisme de $\mathcal{SI}$. En consid\'erant 
$\lambda=\tau^{-1}$, on d\'emontre facilement l'assertion (iii) .

L'assertion (i) est maintenant triviale.\hfill$\Box$

\begin{remark} Une deuxi\`eme bijection naturelle entre $\mathcal{U}$
et $\mathcal{D}$ est donn\'ee par $\alpha\in\mathcal{D}
\longmapsto \alpha'\in\mathcal{U}$. L'application
$$\tau\longmapsto 
\frac{1}{A\circ\left(\int_0A^\tau\right)^{\langle -1\rangle}}$$
(provenant de l'automorphisme $(A,\alpha)\longmapsto
(A(\alpha')^\lambda,\alpha)$ de $\mathcal{SI}$) 
permet d'interpoler entre $\frac{1}{A}$ et la s\'erie r\'eciproque
$$\left(\int_0A\right)^{\langle
  -1\rangle}=\int_0\frac{1}{A\circ\left(\int_0A\right)^{\langle-1\rangle}}$$
de $\int_0A\in\mathcal{SD}$ 
associ\'ee \`a cette deuxi\`eme bijection, voir \cite{Baint}.
\end{remark}

\section{La transform\'ee de Hankel}\label{hankel}

Ce chapitre contient notre r\'esultat principal, sugg\'er\'e
par le th\'eor\`eme \ref{thmA} du chapitre \ref{inverse}.

La $n-$i{\`e}me {\em matrice de Hankel} $H(n)$ d'une suite 
$s=(s_0,s_1,s_2,\,\dots)$ est la matrice sym\'etrique dont les coefficients 
$h_{i,j},\ 0\leq i,j<n$ ne d\'ependent que de la somme $i+j$ des
indices et sont donn\'es par $h_{i,j}=s_{i+j}$. La matrice
$H(n)$ d\'epend donc seulement de $s_0,s_1,\,\dots,s_{2n-2}$.
La {\em transform{\'e}e de Hankel} de $s$ est alors d{\'e}finie
comme {\'e}tant la suite
$$\det(H(1))=s_0,\, \det(H(2))=s_0s_2-s_1^2,\, \det(H(3)),\, \dots\ $$
des d{\'e}terminants des matrices de Hankel d'ordre $1,2,3,\,\dots$ 
associ{\'e}es {\`a} $s$.

Une formule de Hadamard (voir~\cite{G}, page 30, voir aussi \cite{L}) 
implique que deux suites $a$ et $b=I(a)$ dont les s\'eries 
g\'en\'eratrices 
sont reli\'ees par la transformation inverse
$\sum_{n=0}^\infty b_nt^n=\left(\sum_{n=0}^\infty a_n t^n\right)/\left( 1+
t\sum_{n=0}^\infty a_n t^n\right)$, ont m{\^e}me transform{\'e}e 
de Hankel.
Comme les polyn{\^o}mes $I_n(x)$ interpolent les 
it{\'e}r{\'e}es  de la transform{\'e}e 
inverse, la transform{\'e}e de Hankel de la suite
$I^x(a)=(I_0(x),\,I_1(x),\,I_2(x),\,\dots)$
ne d{\'e}pend pas de $x$.

Pour un entier $k>0$, d{\'e}finissons la $k-$i{\`e}me 
transform{\'e}e de Hankel de  $s=(s_0,\,s_1,\,\dots)$ comme
la suite $\big(d_{k,n}=\det(H_k(n))\big)_{n=1,2,\,\dots}=(s_k,s_ks_{k+2}-
s_{k+1}^2,\,\dots)$ o\`u $H_k(n)=(s_{i+j+k})_{0\leq i,j<n}$ est la 
matrice de Hankel de taille $n\times n$ associ\'ee \`a la suite
d\'ecal\'ee $s_k,\,s_{k+1},\,s_{k+2},\,\dots$. 

\begin{thm}\label{laymangen} (i) La suite
$$\Big(\det\big((I_{i+j+k}(x))_{0\leq i,j<n}\big)\Big)_{n=1,2,3,\,\dots}$$
de la $k-$i{\`e}me transform{\'e}e de Hankel de $I^x(a)=(I_0(x),
\,I_1(x),\,\dots)$
ne contient que des polyn{\^o}mes de degr{\'e} $\leq k$ en $x$.

\ \ (ii) Le d\'eterminant 
$$\det((Q_{1+i+j}(x))_{0\leq i,j<n})$$
pour $Q_1,\,Q_2,\,\dots$ associ\'es \`a $s(x)=s_0+s_1x+s_2x^2+\dots$ comme dans 
le chapitre \ref{inverse} ne d\'epend pas de $s_1$.
\end{thm}

\begin{remark}
L'identit\'e de condensation de Dodgson (cf. \cite{Kr}) montre que
les d\'eterminants $d_{k,n}$ v\'erifient l'\'egalit\'e
$$d_{k-1,n+1}\ d_{k+1,n-1}=d_{k-1,n}\ d_{k+1,n}-d_{k,n}^2$$
o\`u l'on a pos\'e $d_{k,0}=1$ pour tout $k$.
Cette identit\'e est parfois utile pour calculer
r\'ecursivement la transform\'ee de Hankel $\left(d_{0,n}\right)_{
n=1,2,\,\dots}$ \`a partir de $d_{k,1}=s_k$.
\end{remark}

\section{Mots de \L ukasiewicz et r\'eversion des s\'eries (Lagrange)}\label{Lukas}

Ce chapitre est d\'evolu \`a l'\'etude des mots de \L ukasiewicz.
Les propri\'et\'es de ces mots sont ensuite exploit\'ees pour 
d\'emontrer les th\'eor\`emes \ref{invlagrange},
\ref{lagrangerev} et \ref{thmA}.

Nous commen\c cons par d\'emontrer 
le th\'eor\`eme \ref{thmA}
qui \'equivaut \`a l'identit\'e
$$[x^k]  \sum_{n=1}^\infty Q_n(x) t^n   =     q(t)^{k+1}\ .$$

Soit
$$
s(x) = s_0 + \sum_{j=1}^\infty s_j x^j
$$
une s\'erie formelle dont les coefficients $s_0, s_1, s_2, \dots$ sont des lettres
qui ne commutent qu'avec la variable $x$.
Comme au d\'ebut, nous associons \`a $s(x)$ la suite des polyn\^omes
$$
P_1(x) =  s_0,
P_2(x) =  s_0^2 + s_0s_1 x,
P_3(x) =  s_0^3 + (s_0^2s_1+s_0s_1s_0)x + (s_0^2s_2+s_0s_1^2)x^2,
\dots
$$
d\'efinie de fa\c con r\'ecursive par $P_k(x) = \lfloor P_{k-1}(x) s(x) \rfloor_k$. Notons $[x^k]P_n(x)$ le coefficient de~$x^k$
du polyn\^ome~$P_n(x)$. On a une bijection entre les mon\^omes de~$[x^k]P_n(x)$
et les chemins sur $\mathbb{N}\times\mathbb{N}$ de $(0,0)$ \`a $(n,k)$
ne traversant pas la diagonale $y=x$ et 
qui n'utilisent
que des pas $(1,0),(0,1)$ orient\'es vers le nord ou vers l'est.
En effet, associons \`a $s_{i_1}s_{i_2}s_{i_3}\dots s_{i_{n}}$ le chemin
$$i_1\times(0,1) + (1,0) + i_2\times(0,1) + (1,0) + i_3\times(0,1) + (1,0) + \cdots + i_{n}\times(0,1) + (1,0)$$
(on a toujours $i_1 = 0$). En particulier, le nombre de tels mon\^omes contribuant au
coefficient $[x^{n-1}]P_n(x)$ de plus haut degr\'e est donn\'e par le
nombre de Catalan $C_{n-1} = \binom{2(n-1)}{n-1} \frac{1}{n}$.
Posons $Q_n(x) = x^{n-1}{\overline P}_n(1/x)$ o\`u ${\overline P}_n(x)$
est obtenu en lisant \`a l'envers les mon\^omes
contribuant aux coefficients $x^0,\,x^1,\,\dots,\,x^{n-1}$ de $P_n(x)$ :
$$
Q_1(x) =  s_0, \quad
Q_2(x) =  s_0^2x + s_1s_0, \quad
Q_3(x) =  s_0^3x^2 + (s_1s_0^2+s_0s_1s_0)x + s_2s_0^2+s_1^2s_0,\,\dots .
$$
Munissons la lettre~$s_i$ du poids $w(s_i) = i-1$ 
et posons
$$w(s_{i_1}s_{i_2}s_{i_3}\dots s_{i_{n}}) =
w(s_{i_1})+w(s_{i_2})+w(s_{i_3})+\dots+w(s_{i_{n}})$$ 
pour un mot $s_{i_1}s_{i_2}s_{i_3}\dots s_{i_{n}}$ de longueur $n$. 
Repr\'esentons un mot $s_{i_1}\dots
s_{i_n}$ apparaissant dans $Q_n(x)$ par le chemin de sommets $$(0,0),\,
(1,s_{i_1}-1),\, (2,s_{i_1}+s_{i_2}-2),\, \dots,\,
(n,-n+\sum_{j=1}^n i_j)=(n,w(s_{i_1}\dots s_{i_n}))$$ obtenu en
concat\'enant les pas $(1,w(s_{i_j}))=(1,{i_j}-1)$ associ\'es \`a
$s_{i_1},\,s_{i_2},\,\dots,\,s_{i_n}$.
\medskip

\centerline{\epsfysize3cm\epsfbox{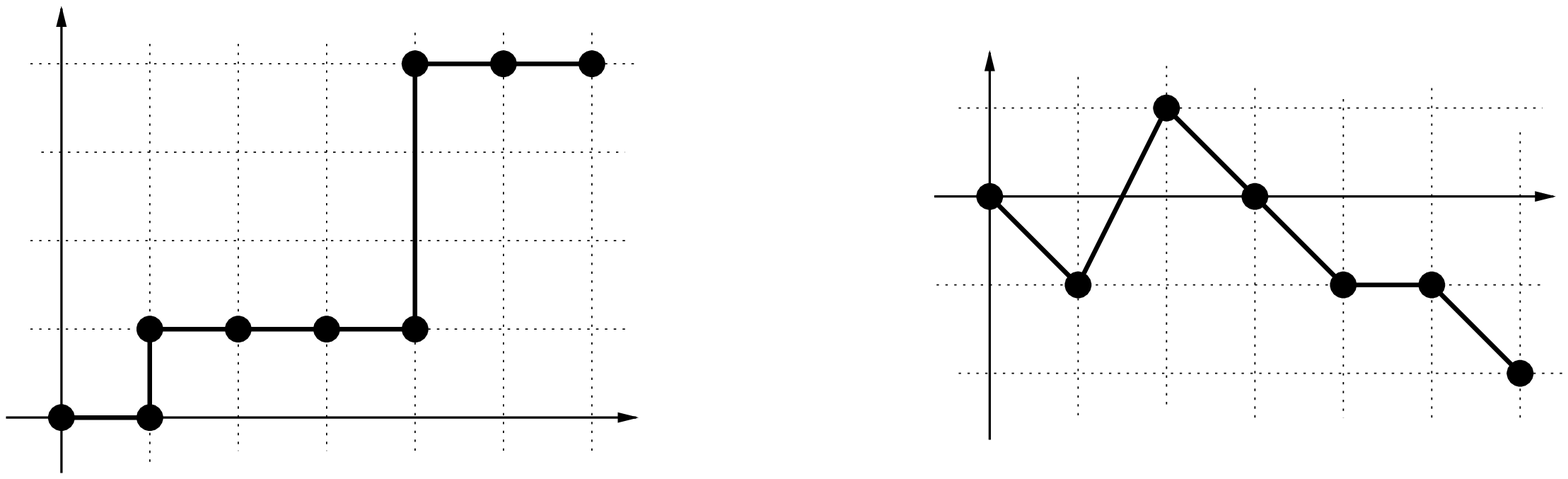}}
\centerline{Le chemin associ\'e au mot $s_0s_1s_0s_0s_3s_0$ (dans $P_6(x)$)
et \`a son miroir $s_0s_3s_0s_0s_1s_0$ (dans $Q_6(x)$).} 
\medskip

Les mots qui apparaissent dans $Q_n(0) = [x^0]Q_n(x)$ sont les
\emph{mots de \L ukasiewicz} (voir \cite{C}). 
Leur s\'erie g\'en\'eratrice est donn\'ee par 
$q(t)=\sum_{n=1}^\infty Q_n(0)t^n$. Remarquons que les mots 
de \L ukasiewicz
de $Q_{n+1}(0)$ sont en bijection avec les parenth\'esages
de longueur $2n$ comportant $n$ parenth\`eses ouvrantes et fermantes.
Pour le voir on commence par supprimer la derni\`ere lettre
$s_0$ d'un mot de \L ukasiewicz et on remplace ensuite
une lettre $s_k$ par le mot de longueur $k+1$
consistant en $k$ parenth\`eses ouvrantes \lq\lq  (($\dots$((''
suivi d'une parenth\`ese fermante \lq\lq  )''. Pour le mot
de \L ukasiewicz $s_2s_0s_1s_2s_2s_0s_0s_0$ on obtient ainsi
$$\begin{array}{cccccccc}
s_2&s_0&s_1&s_2&s_2&s_0&s_0&\\
(()& )& ()& (()& (()& )& )&\end{array}_ .$$  

\begin{lem} \label{caracterisation_des_coeffs}
Le coefficient $[x^k]Q_n(x)$ de $Q_n(x)$ est la 
somme de tous les mots $s_{i_1}s_{i_2}\dots s_{i_{n}}$ 
de longueur $n$ et de poids
$w(s_{i_1}s_{i_2}s_{i_3}\dots s_{i_{n}}) =\sum_{j=1}^n i_j\ -n= -(k+1)$  
tels que $w(s_{i_1}s_{i_2}s_{i_3}\dots s_{i_h}) \ge -k$ pour tout $h <n $.
\end{lem}

\noindent{\bf Preuve.} 
Le coefficient $[x^k]P_n(x)$ est constitu\'e de tous les mots
$s_{i_1}\dots s_{i_n},\ \sum_{j=1}^n i_j=k$, qui v\'erifient les 
in\'egalit\'es $\sum_{j=1}^l i_j\leq l-1$ pour $l=1,\,\dots,n$.

Le miroir $s_{i_n}\dots s_{i_1}$ de poids $w(s_{i_n}\dots s_{i_1})=
\sum_{j=1}^n (i_j-1)=-(n-k)$ d'un tel mot contribue au coefficient
$[x^{n-1-k}]Q_n(x)$. Nous avons 
$$w(s_{i_n}\dots s_{i_{n+1-h}})=-h+\sum_{j=n+1-h}^n i_j=-h+k-
\sum_{j=1}^{n-h}i_j\ .$$
En utilisant la majoration $\sum_{j=1}^{n-h} i_j\leq n-h-1$
rencontr\'ee
ci-dessus, nous avons pour $h<n$
$$w(s_{i_n}\dots s_{i_{n+1-h}})\geq -h+k-(n-1-h)>-(n-k)\ .\hfill \Box$$

\noindent{\bf Preuve du th\'eor\`eme \ref{thmA}.}
Soit $s_{i_1}\dots s_{i_n}$ un mot de longueur $n$ et de poids
$w(s_{i_1} \dots s_{i_n})=\sum_{j=1}^n i_j\ -n=-(k+1)$ 
contri\-buant au coefficient $[x^k]Q_n(x)$ de $Q_n$.
Un tel mot s'\'ecrit de mani\`ere unique sous la forme
$s_{i_1} \dots s_{i_n}=l_0\dots l_k$
o\`u les mots $l_0,\,\dots,\,l_k$ sont des mots de \L ukasiewicz en
$s_0,\,s_1,\,s_2,\,\dots$ (voir la remarque ci-dessous pour un exemple).
En effet, soit $a\geq 1$ le plus petit indice tel que $w(s_{i_1}\dots s_{i_a})=-1$.
Le mot $l_0=s_{i_1}\dots s_{i_a}$ satisfait alors les conditions du lemme
\ref{caracterisation_des_coeffs} avec $k=0$. C'est donc un mot de \L ukasiewicz.
De plus, c'est le seul sous-mot initial de $s_{i_1}\dots s_{i_n}$ qui soit de
\L ukasiewicz car
un sous-mot initial de la forme $s_{i_1}\dots s_{i_b}$ avec $b<a$
est de poids $w(s_{i_1}\dots s_{i_b})\geq 0$. D'autre part, un tel mot
avec $b>a$ ne peut \^etre \`a la fois de poids $-1$ et v\'erifier
les conditions du lemme \ref{caracterisation_des_coeffs}.

Si $k=0$, le lemme \ref{caracterisation_des_coeffs} implique que $a=n$.
Pour $k>0$, on a $a<n$ et le compl\'ement
$s_{i_{a+1}}\dots s_{i_n}$ est un mot de poids $-k$ v\'erifiant de
nouveau les conditions du lemme \ref{caracterisation_des_coeffs}.
Par r\'ecurrence sur $k$, on a alors
$s_{i_{a+1}}\dots s_{i_n}=l_1\dots l_k$ avec $l_1,\,\dots,\,l_k$ des mots de
\L ukasiewicz. Ceci montre que l'ensemble des mots formant le coefficient
$[x^k]Q_n(x)$ est l'ensemble des mots de longueur $n$ en $s_0,\,s_1,\,\dots$
obtenus en concat\'enant $(k+1)$ mots de \L ukasiewicz.
On a donc l'\'egalit\'e
$[x^k]Q_n(x)=[t^n]q(t)^{k+1}$. \hfill $\Box$

\begin{remark} La factorisation $s_{i_1}\dots s_{i_n}=l_0\dots l_k$
d'un mot de poids $-(k+1)$ satisfaisant les conditions du lemme
\ref{caracterisation_des_coeffs} en
$(k+1)$ mots de \L ukasiewicz est bien visible sur la repr\'esentation
graphique introduite ci-dessus. Ainsi, pour le mot $s_0s_3s_0s_0s_1s_0$
contribuant au coefficient $x^{5-(0+3+0+0+1+0)}=x$ de $Q_6(x)$, on obtient
$l_0=s_0$ et $l_1=s_3s_0s_0s_1s_0$.
\end{remark}

\noindent {\bf Preuve du th\'eor\`eme \ref{invlagrange}.}
Soit  $s_{i_1}s_{i_2}s_{i_3}\dots s_{i_{n}}$ un mot de  \L ukasiewicz.
Si $n=1$, alors le mot est \'egal \`a $s_0$. Si $n \ge 2$, alors $i_1 \ge 1$
et $s_{i_2}s_{i_3}\dots s_{i_{n}}$ est un mon\^ome de~$[x^{i_1-1}] Q_{n-1}(x)$.
Il se factorise donc en $i_1$ facteurs  de \L ukasiewicz.

Ceci sugg\`ere de consid\'erer la bijection suivante entre les  mots
de  \L ukasiewicz et les arbres plans
enracin\'es: Au mot de  \L ukasiewicz $s_{i_1}s_{i_2}s_{i_3}\dots s_{i_{n}}$
on fait correspondre l'arbre plan de $n$~sommets muni d'une racine de
degr\'e~$i_1$. Les $i_1$ fils de la racine correspondent r\'ecursivement
aux $i_1$ facteurs de \L ukasiewicz du mot $s_{i_2}s_{i_3}\dots s_{i_{n}}$.
Cette bijection se traduit par l'identit\'e
$
q(t) = ts(q(t))
$
pour les s\'eries g\'en\'eratrices.
\hfill$\Box$

\medskip

\centerline{\epsfysize4cm\epsfbox{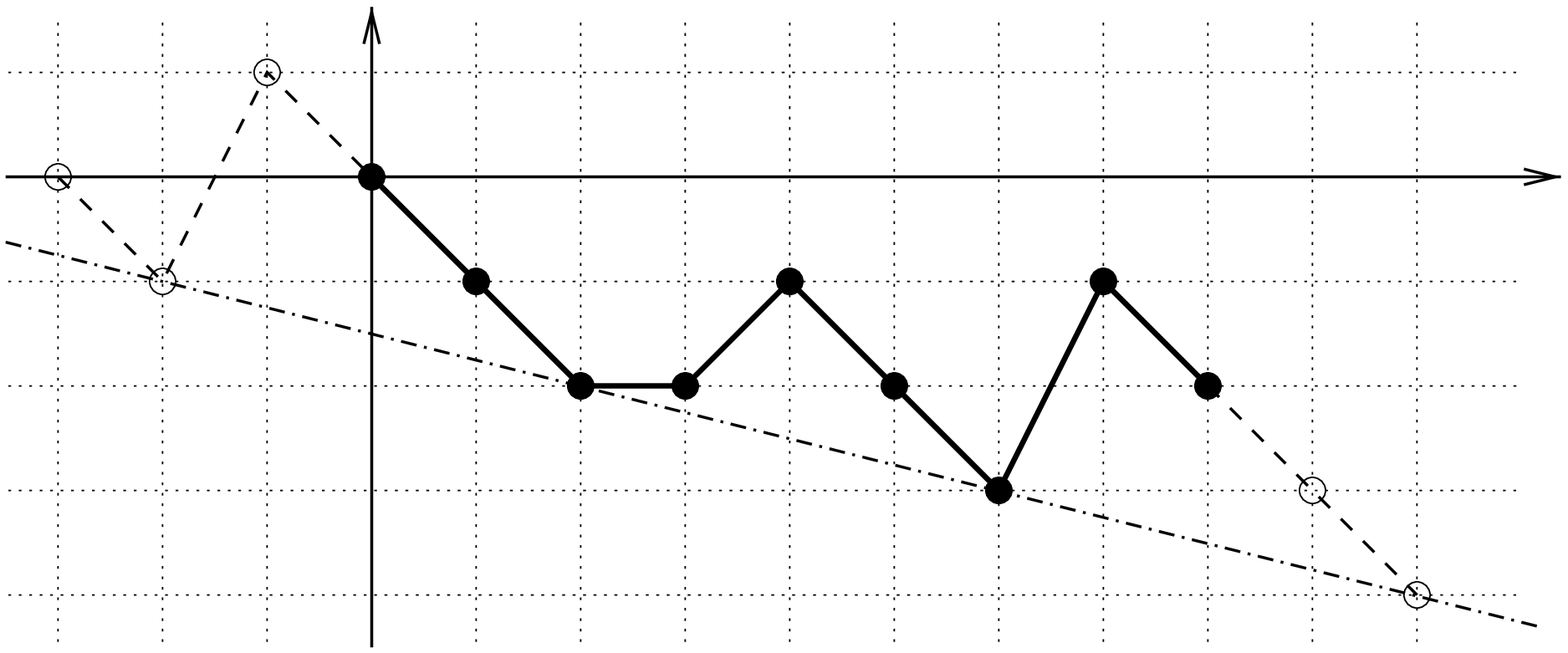}}
\centerline{Le mot
$s_0s_0s_1s_2s_0s_0s_3s_0$, rendu cyclique.} 
\medskip

\noindent{\bf Preuve du th\'eor\`eme \ref{lagrangerev}.}
Le coefficient $[x^{n-k-1}]s(x)^n$ compte tous les mots
$s_{i_1}s_{i_2}\dots s_{i_n}$ de longueur $n$ qui sont de poids
$w(s_{i_1}s_{i_2}\dots s_{i_n})=\sum_{j=1}^ni_j\ -n=(n-k-1)-n=-(k+1)$.
Ainsi le mot $s_0s_0s_1s_2s_0s_0s_3s_0$ apparaissant dans
$[x^6]s(x)^8$ est de longueur $8$ et de poids
$-2$ et illustre le cas particulier $n=8$ et $k=1$.
Associons \`a un tel mot la suite $n-$p\'eriodique (bi-infinie) 
de lettres
$\dots s_{i_{-2}}s_{i_{-1}}s_{i_0}
s_{i_1}s_{i_2}\dots s_{i_n}s_{i_{n+1}}s_{i_{n+2}}\dots$
avec $s_{i_{n+h}} = s_{i_h}$ pour tout $h \in {\mathbf Z}$. 
Regardons la
repr\'esentation graphique, c'est-\`a-dire la suite infinie de points 
$$\begin{array}{l}\displaystyle
\dots,\, (-2,-w(s_{i_{-1}}s_{i_0})),\, (-1,-w(s_{i_0})),\, (0,0),\,
(1,w(s_{i_1})),\, (2,w(s_{i_1}s_{i_2})),\  \dots,\\
\displaystyle
(n, w(s_{i_1}s_{i_2}\dots s_{i_n})) = (n,-(k+1)),\ 
(n+1, w(s_{i_1}s_{i_2}\dots s_{i_{n+1}})) = (n+1,-(k+1)+w(s_{i_1})),\\
\displaystyle
(n+2, w(s_{i_1}s_{i_2}\dots s_{i_{n+2}})) = (n+2,-(k+1)+w(s_{i_1}s_{i_2})),\,  \dots\end{array}$$
La suite des produits scalaires de ces points avec le vecteur $(k+1,n)$ est p\'eriodique
et la longueur de la p\'eriode est un diviseur de~$n$. Supposons que la valeur minimale
de ces produits scalaires est prise sur le point $(h, w(s_{i_1}s_{i_2}\dots s_{i_h}))$
avec $h \in \{1,2,\,\dots,\,n\}$. Par ailleurs, $h$ est unique 
(modulo $n$) si $k = 0$. De toute
fa\c con, le mot $s_{i_{h+1}}s_{i_{h+2}}\dots s_{i_n}s_{i_1}s_{i_2}\dots s_{i_h}$
appara\^\i t dans $[x^k]Q_n(x)$. On appelle ce mot un \emph{r\'earrangement cyclique}
de $s_{i_1}s_{i_2}\dots s_{i_n}$. D'apr\`es la d\'emonstration du th\'eor\`eme 1,
ce r\'earrangement cyclique a une factorisation canonique en $k+1$ mots de
\L ukasiewicz. Parmi les $n$~r\'earrangements  cycliques possibles de
$s_{i_1}s_{i_2}\dots s_{i_n}$ il y en a donc exactement
$(k+1)$ qui apparaissent dans $[x^k]Q_n(x)$~: les $(k+1)$~r\'earrangements
cycliques des  facteurs de \L ukasiewicz de
$s_{i_{h+1}}s_{i_{h+2}}\dots s_{i_n}s_{i_1}s_{i_2}\dots s_{i_h}$.
En effet, si l'on choisit un r\'earrangement  cyclique dont la premi\`ere lettre
n'est pas la premi\`ere lettre d'un facteur de \L ukasiewicz, alors l'in\'egalit\'e
n\'ecessaire pour l'appartenance \`a  $[x^k]Q_n(x)$ n'est pas satisfaite pour le
mot qui va jusqu'\`a la derni\`ere lettre du facteur pr\'ec\'edent.

Ainsi, pour notre exemple $s_0s_0s_1s_2s_0s_0s_3s_0$ repr\'esent\'e par
la figure ci-dessus, les facteurs de 
\L ukasiewicz du mot cyclique bi-infini sont d\'elimit\'es par les 
intersections du graphe repr\'esentant ce mot avec la droite $4y=-6-x$,
repr\'esent\'ee en pointill\'e. Ses deux facteurs de \L ukasiewicz sont donc
$s_1s_2s_0s_0$ et $s_3s_0s_0s_0$. Parmi les huit r\'earrangement 
circulaire du mot $s_0s_0s_1s_2s_0s_0s_3s_0$, il n'y a donc que 
$s_1s_2s_0s_0s_3s_0s_0s_0$ et $s_3s_0s_0s_0s_1s_2s_0s_0$ qui 
apparaissent dans $[x^1]Q_8(x)$.

Dans le cas g\'en\'eral, on obtient ainsi une bijection
\begin{eqnarray*}
\{1,\,\dots,\, k+1\}\times \left\{
\hbox{\begin{tabular}{p{3cm}}
mots en $s_0,\,s_1,\,\dots$ de lon\-gueur $n$ et de poids $-(k+1)$
\end{tabular}} \right\}
&\to&
\{1,\dots ,n\}\times \left\{
\hbox{\begin{tabular}{p{4.5cm}}
produits de $(k+1)$ mots de \L u\-ka\-sie\-wicz de lon\-gueur to\-tale $n$ en $s_0,s_1,\dots$
\end{tabular}}\right\} \\
(k',mot) &\mapsto& (n',\,luk),
\end{eqnarray*}
o\`u $luk$ est le r\'earrangement cyclique de $mot$ qui appartient \`a $[x^k]Q_n(x)$
et qui fait appara\^\i tre la premi\`ere lettre de $mot$ dans le $k'$-ième facteur
de \L ukasiewicz de $luk$ ($n'$ correspond \`a la nouvelle place de la
premi\`ere lettre de $mot$). Cette bijection implique l'\'egalit\'e
$$(k+1)\ [x^{n-k-1}] s(x)^n=n\ [t^n]q(t)^{k+1}.\qquad \Box$$

\begin{remark} Dans le contexte d'une variable $t$ ne commutant 
pas avec les variables $s_i$, il faudrait introduire la
variable~$t$ devant chaque lettre, i.e.
\begin{eqnarray*}
q(t) &=&    ts_0 + ts_1 ts_0 + ts_2 (ts_0)^2 + (ts_1)^2 ts_0 +  \\
        &  & {} ts_3 (ts_0)^3 +  ts_2 ts_1 (ts_0)^2 +  ts_2 ts_0 ts_1 ts_0 +   ts_1 ts_2   (ts_0)^2 + (ts_1)^3 ts_0
                     +\cdots
\end{eqnarray*}
\end{remark}

La derni\`ere lettre d'un  mot de \L ukasiewicz  $s_{i_1}s_{i_2}s_{i_3}\dots s_{i_{n}}$
est toujours la lettre $s_0$  (i.e. $i_{n} = 0$).  Comme  $w(s_1) =
0$, le mot en les lettres $s_0,\,s_2,\,s_3,\,\dots$ obtenu 
par suppression de toutes les lettres $s_1$ dans un mot
de \L ukasiewicz est encore  un  mot
de \L ukasiewicz. Appelons-le mot \emph{r\'eduit}  de \L ukasiewicz
et notons $q_{s_1=0} (t)$ la s\'erie g\'en\'eratrice des mots
r\'eduits de \L ukasiewicz. Nous avons alors le r\'esultat suivant, 
utile au chapitre \ref{sectdethankel}.

\begin{prop} \label{indeps1}
On a l'\'egalit\'e entre s\'eries g\'en\'eratrices
$\;\;  q(t) \;=\;    q_{s_1=0} ((1-ts_1)^{-1}\ t)$.
\end{prop}

\noindent{\bf Preuve.}
Un  mot de \L ukasiewicz  $l=s_{i_1}s_{i_2}s_{i_3}\dots s_{i_{n}}$ qui
ne contient pas la
lettre~$s_1$, est le mot r\'eduit  de \L ukasiewicz pour tous
les mots de  \L ukasiewicz de la forme
$s_1^{k_1}  s_{i_1} s_1^{k_2} s_{i_2} s_1^{k_3}  s_{i_3} \dots s_1^{k_{n}}   s_{i_{n}}$
avec $k_1,\,k_2,\,k_3,\,\dots,\,k_{n} \in \{0,1,2,\,\dots\}$.

Le mot r\'eduit $l$ intervient avec une contribution de
$ts_{i_1}ts_{i_2}ts_{i_3}\dots ts_{i_{n}}$ dans
la s\'erie g\'en\'eratrice $q_{s_1=0}(t)$ des mots r\'eduits.
L'ensemble
de tous les mots de \L ukasiewicz dont $l$ est le mot r\'eduit
contribue donc avec $(1-ts_1)^{-1}ts_{i_1}(1-ts_1)^{-1}ts_{i_2}\dots
(1-ts_1)^{-1}ts_{i_{n}}$ \`a la
s\'erie g\'en\'eratrice $q(t)$ de tous les mots de \L ukasiewicz.\hfill$\Box$

\subsection{Digression: Arbres binaires r\'eguliers,
arbres plans enracin\'es et mots de \L ukasiewicz} Un {\it arbre 
binaire r\'egulier} est un arbre plan enracin\'e 
(modulo la relation d'\'equivalence \'evidente) dont tous les 
sommets ont z\'ero ou deux enfants. 
Notons ${\mathcal B}_n$ l'ensemble des arbres binaires 
r\'eguliers avec $n+1$ feuilles (et $2n+1$ sommets, $2n$ ar\^etes)
et ${\mathcal T}_n$ l'ensemble
des arbres plans enracin\'es ayant $n+1$ sommets (et $n$ ar\^etes).
Les deux ensembles ${\mathcal B}_n,\ {\mathcal T}_n$ ont m\^eme cardinalit\'e,
donn\'ee par le $n-$i\`eme nombre de Catalan ${2n\choose n}/(n+1)$
(cf. eg. l'Exercice 6.19 d,e dans \cite{Stanl}). Une bijection entre ces deux ensembles finis peut
\^etre d\'ecrite comme suit: Un arbre binaire r\'egulier
$B\in{\mathcal B}_n$ poss\`ede exactement $n$ ar\^etes gauches
(orient\'ees NO) et $n$ ar\^etes droites (orient\'ees NE). En contractant
toutes les ar\^etes gauches (respectivement droites) de $B$ on obtient
un arbre planaire enracin\'e $C_L(B)$ (respectivement $C_R(B)$) dans
${\mathcal T}_n$ et on v\'erifie facilement que les deux applications
$C_L,C_R:{\mathcal B}_n\longrightarrow {\mathcal T}_n$ sont bijectives. 

\medskip

\centerline{\epsfysize5cm\epsfbox{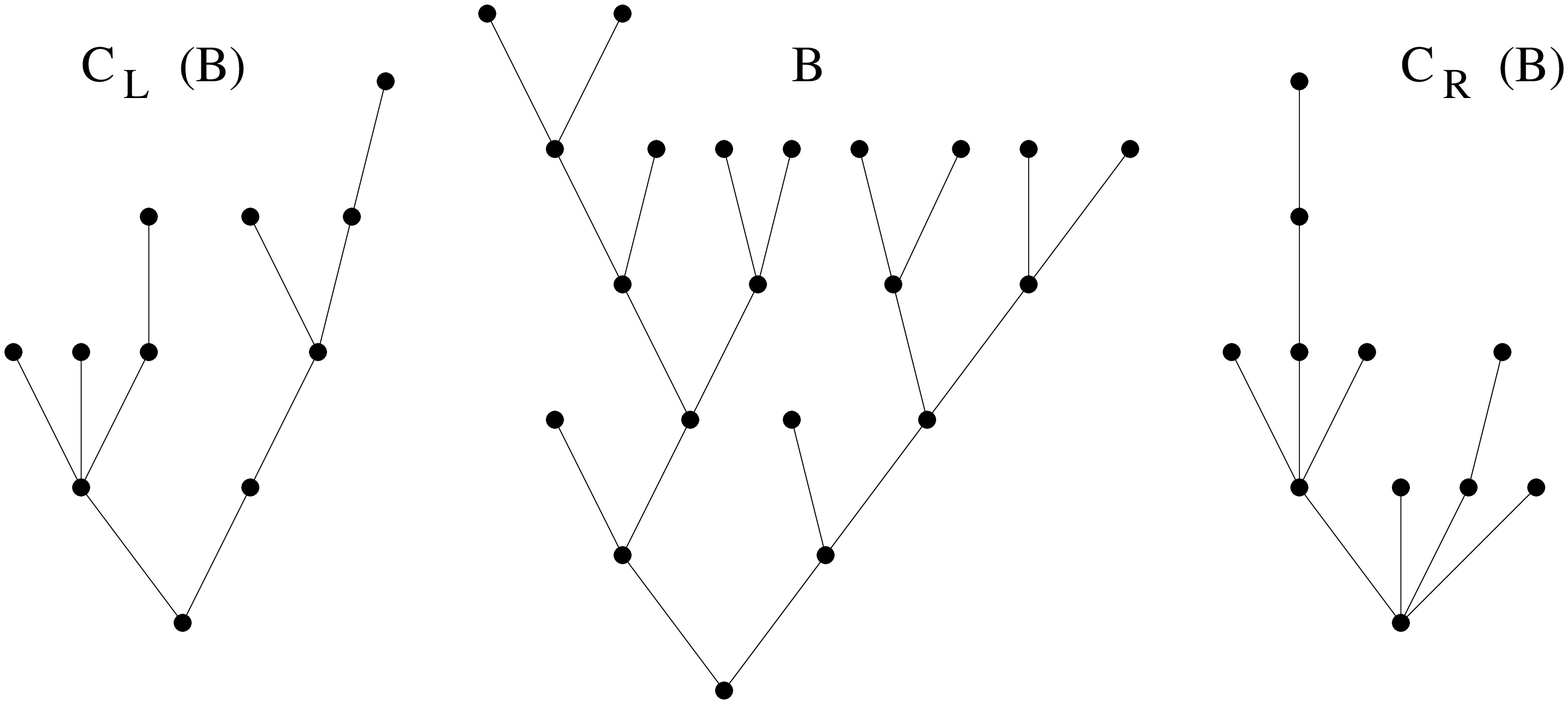}}
\centerline{Les deux arbres $C_L(B),C_R(B)$
associ\'es \`a un arbre binaire $B$.} 
\medskip

D\'esignons par $\overline X$ l'arbre 
``miroir'' obtenu en r\'eflechissant un arbre  
$X\in {\mathcal B}_n$ ou $X\in {\mathcal T}_n$ par rapport 
\`a une droite verticale. On montre facilement l'identit\'e 
$\overline{C_R(B)}=C_L(\overline B)$. En conjugant l'involution 
$T\longmapsto \overline T$ sur ${\mathcal T}_n$ par les
bijections $C_R,C_L$, on obtient ainsi deux involutions 
$\iota_R(B)=C_R^{-1}(\overline{C_R(B)})=\overline{C_L^{-1}(C_R(B))}$ et
$\iota_L(B)=C_L^{-1}(\overline{C_L(B)})=\overline{C_R^{-1}(C_L(B))}$ sur 
${\mathcal B_n}$. Une construction analgoue, \`a savoir $\tilde \iota_R=
C_R(\overline{C_R^{-1}(T)})$ et $\tilde \iota_L=
C_L(\overline{C_L^{-1}(T)})$ d\'efinit deux involutions $\tilde \iota_R,
\tilde \iota_L$ sur ${\mathcal T}_n$.
Il serait int\'eressant de comprendre les orbites dans 
${\mathcal B}_n$ (respectivement ${\mathcal T}_n$)
sous l'action du groupe di\'edral de g\'en\'erateurs
$\iota_R,\iota_L$ (respectivement $\tilde\iota_R,\tilde \iota_L$).
En particulier, les points fixes de $\iota_R$ (ou les points fixes
de $\iota_L$) sont en bijection avec les arbres ``sym\'etriques'' de 
${\mathcal T}_n$ qui satisfont $\overline T=T$ (au nombre de
${n\choose \lfloor n/2\rfloor}$) tandis que les points fixes de 
$\tilde \iota_R$  (ou les points fixes
de $\tilde \iota_L$) correspondent bijectivement aux arbres
sym\'etriques binaires de ${\mathcal B}_n$. Le nombre d'arbres 
sym\'etriques binaires r\'eguliers est donn\'e par le nombre de
Catalan ${2m\choose m}/(m+1)$ pour
$n=2m+1$ impair. Pour $n>0$ pair de tels arbres n'existent pas.

Pour terminer cette digression, mentionnons encore le fait (d\'ej\`a 
rencontr\'e dans la preuve du th\'eor\`eme \ref{invlagrange}) que la suite
$i_1,\,\dots,\,i_{n+1}$
des valences des $n+1$ sommets rencontr\'es pour la premi\`ere 
fois lorsqu'on 
contourne un arbre $T\in {\mathcal T}_n$ en partant de sa racine 
d\'efinit bijectivement un mot de \L ukasiewicz $s_{i_1}\dots s_{i_{n+1}}$ 
de longueur $n+1$.

\section{D\'eterminants de Hankel}\label{sectdethankel}

Le but de ce chapitre est la preuve du th\'eor\`eme \ref{laymangen}.
Pour cela, nous introduisons les mots de Motzkin et rappelons 
quelques-unes  de leurs propri\'et\'es. Des \'etudes plus
compl\`etes sont contenues par exemple dans 
\cite{F} et \cite{V}, voir aussi \cite{V1}.

Un {\it chemin de Motzkin} de longueur $n$ est un chemin dans le
premier quadrant $x,y\geq 0$ qui relie l'origine $(0,0)$ au point
$(n,0)$ en utilisant $n$ pas de la forme $(1,-1),\ (1,0)$ ou $(1,1)$.

\medskip
\centerline{\epsfysize2.8cm\epsfbox{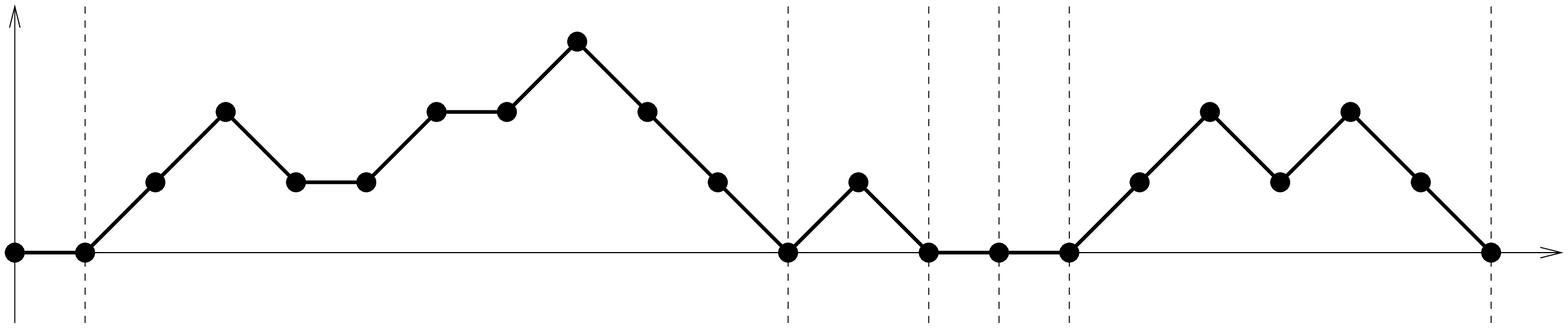}}
\centerline{Un chemin 
de Motzkin de longueur $21$ et ses $6$ facteurs premiers.} 
\medskip

Notons $\Gamma(n)$
l'ensemble des chemins de Motzkin de longueur~$n$. Chaque $\gamma \in \Gamma(n)$ est
affect\'e d'un poids $w(\gamma)$ d\'efini comme le produit des poids des diff\'erents
arcs qui le constituent~: un palier~$((i,h),(i+1,h))$ situ\'e \`a la hauteur~$h$ est
affect\'e du poids~$p(h)$; une descente~$((i,h+1),(i+1,h))$  de la hauteur~$h+1$
\`a la hauteur~$h$ est affect\'e du poids~$q(h)$; enfin, chaque
mont\'ee~$((i,h),(i+1,h+1))$ est affect\'e du poids~$1$.
Le poids~$w(\gamma)$ est ainsi un mon\^ome en les variables (commutatives)
 $p(0)$, $p(1)$, $p(2)$, \dots et $q(0),\ q(1),\,q(2),\,\dots$ et on peut
former la s\'erie g\'en\'eratrice
$$\begin{array}{lcl}
\displaystyle c(u)&=&
1 + \sum_{n=1}^\infty u^n \sum_{\gamma \in \Gamma(n)} w(\gamma)\\
\displaystyle &=&
1 + p(0)u + (p(0)^2 +q(0))u^2 + (p(0)^3+2p(0)q(0)+p(1)q(0))u^3  + \cdots
\end{array}
$$
des chemins de Motzkin. Le chemin de Motzkin de longueur 
$21$ repr\'esent\'e ci-dessus contribue ainsi avec
$$\begin{array}{l}
\displaystyle p(0)\ 1\ 1\ q(1)\ p(1)\ 1\ p(2)\ 1\ q(2)\ q(1)\ q(0)\ 1\ q(0)\ p(0)\ 
p(0)\ 1\ 1\ q(1)\ 1\ q(1)\ q(0)\\
\displaystyle \qquad =p(0)^3p(1)p(2)q(0)^3q(1)^4q(2)\end{array}$$
au coefficient $[u^{21}]c(u)$.

Un  chemin de Motzkin est \emph{premier} s'il n'intersecte la droite
horizontale discr\`ete ${\mathbb Z}\times \{0\}$ qu'en 
ses extr\'emit\'es $(0,0)$ et $(n,0)$. Il est clair que tout
chemin de Motzkin premier est soit un palier \`a la hauteur~$0$ (valu\'e~$p(0)$),
soit commence avec une mont\'ee de la hauteur~$0$ \`a la hauteur~$1$, continue
avec un chemin de Motzkin (\'eventuellement vide) allant de la hauteur~$1$ \`a
la hauteur~$1$ et se termine avec une descente  de la hauteur~$1$ \`a la hauteur~$0$
(valu\'ee~$q(0)$). De plus, tout chemin de Motzkin non vide se
factorise de mani\`ere unique en produit de chemins de Motzkin premiers
(il suffit de consid\'erer les sommets situ\'es \`a la hauteur~$0$).
En it\'erant, on obtient imm\'ediatement le th\'eor\`eme suivant
(voir~\cite{F}).

\begin{thm} Soit
$$
c(u) \; = \; 1 + \sum_{n=1}^\infty c_n u^n
\; = \;
1 + \sum_{n=1}^\infty u^n \sum_{\gamma \in \Gamma(n)} w(\gamma)
$$
la fonction g\'en\'eratrice des chemins de Motzkin. Alors
$$
c(u)
\; = \;
\frac{1}{1-p(0)u-\displaystyle\frac{q(0)u^2}{1-p(1)u-\displaystyle\frac{q(1)u^2}{1-p(2)u-
\displaystyle\frac{q(2)u^2}{\ddots}}}}
$$
\end{thm}

On appelle le d\'eveloppement du th\'eor\`eme pr\'ec\'edent {\it fraction continue de Jacobi},
ou encore \hbox{\emph{J-fraction}}. Il permet d'exprimer les coefficients~$c_n$  d'une s\'erie formelle
\`a l'aide de chemins (de Motzkin). En fait, on a la g\'en\'eralisation suivante. Soit
$$
d(u) \; = \; d_0 + \sum_{n=1}^\infty d_n u^n
\; = \;
\frac{d_0}{1-p(0)u-\displaystyle\frac{q(0)u^2}{1-p(1)u-\displaystyle\frac{q(1)u^2}{\ddots}}}\; =\; d_0\ c(u)
$$
et soit $D = (d_{i+j})_{0\le i,j < \infty}$ la matrice de Hankel
(infinie) associ\'ee
\`a la suite $d_0,\,d_1,\,d_2,\,\dots$ de s\'erie g\'en\'eratrice~$d(u)$.
Appelons un mineur de~$D$ un \emph{d\'eterminant de Hankel}. 
Un tel mineur sera not\'e
$D\binom{\alpha_0,\,\alpha_1,\,\dots,\,\alpha_k}{\beta_0,\,\beta_1,
\,\dots,\,\beta_k}$,
en d\'esignant par $0 \le \alpha_0 < \alpha_1 < \dots < \alpha_k$ et
$0 \le \beta_0 < \beta_1 < \dots < \beta_k$ les indices respectifs des
lignes et colonnes du mineur extrait. Le coefficient $m_{i,j}$ de
la sous-matrice associ\'ee \`a  
$D\binom{\alpha_0,\,\alpha_1,\,\dots,\,\alpha_k}{\beta_0,\,\beta_1,
\,\dots,\,\beta_k}$
est donc donn\'e par $m_{i,j}=d_{\alpha_i+\beta_j}$ pour $0\leq i,j\leq k$.
Regardons, pour   $0 \le  i \le k$, les points $A_i = (-\alpha_i,0)$ et $B_i = (\beta_i,0)$.
La somme des valuations (relativement aux variables   $p(0),\,p(1),
\,p(2), \dots$ et $q(0),\,q(1),\,q(2),\ \dots$) des chemins de Motzkin allant de $A_i$ \`a $B_j$ est
$c_{\alpha_i+\beta_j}$, le terme $(i,j)$ du d\'eterminant
$C\binom{\alpha_0,\,\alpha_1,\,\dots,\,\alpha_k}{\beta_0,\,\beta_1,\,
\dots,\,\beta_k}$,
o\`u $C = (c_{i+j})_{0\le i,j < \infty}$ est la matrice de Hankel associ\'ee
\`a la s\'erie g\'en\'eratrice $c(u)$. 
On peut donc \'enoncer le th\'eor\`eme suivant
(voir~\cite{V}, chapitres IV et V, \cite{V1} ou \cite{FZ}).

\begin{thm} \label{detchemins}
On a
$$\begin{array}{rcl}
\displaystyle
D\binom{\alpha_0,\,\alpha_1,\,\dots,\,\alpha_k}{\beta_0,\,\beta_1,
\,\dots,\,\beta_k}
&=&
\displaystyle
d_0^{k+1}  C\binom{\alpha_0,\,\alpha_1,\,\dots,\,\alpha_k}{\beta_0,
\,\beta_1,\,\dots,\,\beta_k}\\
\displaystyle &=&\displaystyle
d_0^{k+1}\sum_{(\sigma;\,\gamma_0,\,\gamma_1,\,\dots,\,\gamma_k)}(-1)^{\inv(\sigma)}
w(\gamma_0)w(\gamma_1)\cdots w(\gamma_k),
\end{array}
$$
o\`u  la sommation est \'etendue aux paires form\'ees par
une permutation $\sigma\in S_{k+1}$ et une confi\-guration 
$(\gamma_0,\,\gamma_1,\,\dots,\,\gamma_k)$ de $k+1$~chemins de Motzkin
sans sommets communs avec $\gamma_i$ reliant $A_i$
\`a $B_{\sigma(i)}$ pour tout $0 \le i \le k$.
\end{thm}

\begin{remark} Dans le th\'eor\`eme pr\'ec\'edent,
deux chemins de Motzkin $\gamma_i$ et $\gamma_j$ peuvent 
s'inter\-secter en des points de la forme 
${\mathbf Z}^2+(\frac{1}{2},\frac{1}{2})$ (de tels points ne 
sont pas consid\'er\'es comme \'etant des sommets).
\end{remark}

\noindent{\bf Preuve.} En omettant la condition
\lq\lq sans sommets communs'', on voit que la somme est
$C\binom{\alpha_0,\,\alpha_1,\,\dots,\,\alpha_k}{\beta_0,\,\beta_1,
\,\dots,\,\beta_k}$
par d\'efinition du d\'eterminant. Cependant, si deux chemins
$\gamma_i\not=\gamma_j$ 
ont un sommet commun, alors on peut continuer, \`a partir du premier
sommet commun rencontr\'e, le premier chemin sur le
second et le second sur le premier. Il est \'evident que les contributions de ces deux configurations
s'annulent.\hfill $\Box$

Le th\'eor\`eme pr\'ec\'edent permet de calculer le d\'eterminant
de certaines matrices. On d\'enombre pour cela 
des chemins de Motzkin convenablement pond\'er\'es.

\begin{thm} \label{degp0} On a
$$
\deg_{p(0)} D\binom{\alpha_0,\,\alpha_1,\,\dots,\,\alpha_k}{\beta_0,
\,\beta_1,\,\dots,\,\beta_k}
\; = \;
\deg_{p(0)}  C\binom{\alpha_0,\,\alpha_1,\,\dots,\,\alpha_k}{\beta_0,
\,\beta_1,\,\dots,\,\beta_k}
\; \le \; (\alpha_k -k) + (\beta_k -k)
$$
et
\begin{eqnarray*}
[p(0)^{\alpha_k+\beta_k -2k}]
D\binom{\alpha_0,\,\alpha_1,\,\dots,\,\alpha_k}{\beta_0,\,\beta_1,
\,\dots,\,\beta_k}
&=&
d_0^{k+1}  [p(0)^{\alpha_k+\beta_k -2k}]
C\binom{\alpha_0,\,\alpha_1,\,\dots,\,\alpha_k}{\beta_0,\,\beta_1,
\,\dots,\,\beta_k} \\
&=&
d_0^{k+1}
B\binom{\alpha_0,\,\alpha_1,\,\dots,\,\alpha_{k-1}}{\beta_0,\,\beta_1,
\,\dots,\,\beta_{k-1}}
\end{eqnarray*}
avec  $B = (b_{i+j})_{0\le i,j < \infty}$  et
$$
b_0 + \sum_{n=1}^\infty b_n u^n
\; = \;
\frac{q(0)}{1-p(1)u-\displaystyle\frac{q(1)u^2}{1-p(2)u-\displaystyle\frac{q(2)u^2}{\ddots}}}
$$
En particulier,
$$
D\binom{0,1,\,\dots,\,k}{0,1,\,\dots,\,k}
\; = \;
d_0^{k+1} q(0)^k q(1)^{k-1} q(2)^{k-2} \cdots q(k-2)^2 q(k-1)
$$
ne d\'epend pas de $p(0),\,p(1),\,p(2),\, \dots$
\end{thm}

\noindent{\bf Preuve.} Pour que le degr\'e
$\deg_{p(0)} D\binom{\alpha_0,\,\alpha_1,\,\dots,\,\alpha_k}{\beta_0,
\,\beta_1,\,\dots,\,\beta_k} =
\deg_{p(0)}  C\binom{\alpha_0,\,\alpha_1,\,\dots,\,\alpha_k}{\beta_0,
\,\beta_1,\,\dots,\,\beta_k}$
en la variable~$p(0)$  devienne maximal, il faut (et il suffit) que les configurations
$(\gamma_0,\,\gamma_1,\,\dots,\,\gamma_k)$ satisfassent aux conditions suivantes~:
Le chemin $\gamma_0$ relie $A_0$ \`a~$B_0$ et reste toujours
\`a la hauteur~$0$ (ceci donne une contribution de~$p(0)^{\beta_0+\alpha_0}$).
Le chemin $\gamma_i$ commen\c cant 
\`a~$A_i = (-\alpha_i,0)$, $i \in \{1,2,\,\dots,\,k\}$,
reste \`a la hauteur~$0$  jusqu'\`a $(-\alpha_{i-1}-1,0)$ et monte \`a
$(-\alpha_{i-1},1)$ ensuite (ceci donne une contribution
de~$p(0)^{\alpha_i - \alpha_{i-1} -1}$).
Similairement, le chemin qui se termine 
\`a~$B_i = (\beta_i,0)$,  $i \in \{1,2,\,\dots,\,k\}$,
descend de $(\beta_{i-1},1)$ vers $(\beta_{i-1}+1,0)$
et reste ensuite jusqu'\`a~$B_i$ \`a la hauteur~$0$
(ceci donne une contribution
de~$q(0) p(0)^{\beta_i - \beta_{i-1} -1}$).
Ce qui n'a pas encore \'et\'e consid\'er\'e n'est rien d'autre qu'une  confi\-guration
de $k$~chemins de Motzkin (translat\'es par le vecteur $(0,1)$)
deux \`a deux disjoints qui relient les sommets
$(-\alpha_0,1),\,(-\alpha_1,1),\,  \dots,\,(-\alpha_{k-1},1)$
aux sommets  $(\beta_0,1),\,(\beta_1,1),\,\dots,\,(\beta_{k-1},1)$.
\hfill$\Box$
\medskip

\noindent{\bf Preuve de l'assertion (i) du th\'eor\`eme \ref{laymangen}.}
Soit $a(u) = a_0 + a_1 u + a_2 u^2 + a_3 u^3 +a_4 u^4 +\cdots$ une 
s\'erie g\'en\'eratrice. Comme les termes de la suite 
$\Big(\det\big((I_{i+j+k}(x))_{0\leq i,j<n}\big)\Big)_{n=1,2,3,\,\dots}$
d\'ependent polynomialement de $x$ et des coefficients $a_0,\,a_1,\,\dots$,
il suffit d'\'etudier le cas g\'en\'erique $a_0\not= 0$.
Consid\'erons la transform\'ee inverse continue
$I^x[a(u)] =  a(u)/\bigl(1+xua(u)\bigr)$.  On a alors
$$
I^x[a(u)] \; = \;  \frac{a_0}{\frac{a_0}{a(u)} + a_0 u x} \; = \;
\frac{a_0}{1 - \left(\frac{a_1}{a_0} - a_0 x\right) u -
\left( \frac{1}{u^2} - \frac{a_1}{a_0 u} - \frac{a_0}{a(u) u^2}\right)u^2}
$$
et on remarque que
$$ \frac{1}{u^2} - \frac{a_1}{a_0 u} - \frac{a_0}{a(u) u^2}=
\frac{a_0a_2-a_1^2}{a_0^2}+\frac{a_0^2a_3-2a_0a_1a_2+a_1^3}{a_0^3}u+
\dots$$
n'a pas de p\^ole en $0$. On cherche \`a calculer le degr\'e en $x$
de $D\binom{0,1,\,\dots,\,n-1}{k,\,k+1,\,\dots,\,k+n-1}$ o\`u $D$
est comme ci-dessus. Le th\'eor\`eme
\ref{degp0} montre que ce degr\'e est $\leq (n-1)-(n-1)+k+n-1-(n-1)=k$
ce qui prouve l'assertion (i) du th\'eor\`eme \ref{laymangen}.
\hfill $\Box$

\begin{remark} Soit $A_x$ la matrice de Hankel associ\'ee \`a la 
s\'erie~$I^x[a(u)]$. Le th\'eor\`eme \ref{degp0} montre
l'\'egalit\'e
$$
[x^{\alpha_k+\beta_k -2k}]
A_x \binom{\alpha_0,\,\alpha_1,\,\dots,\,\alpha_k}{\beta_0,\,\beta_1,
\,\dots,\,\beta_k}
\; = \;
a_0^{k+1}  (-a_0)^{\alpha_k+\beta_k -2k}
\widehat{A}
\binom{\alpha_0,\,\alpha_1,\,\dots,\,\alpha_{k-1}}{\beta_0,
\,\beta_1,\,\dots,\,\beta_{k-1}}
$$
o\`u $\widehat{A}$ est   la matrice de Hankel associ\'ee \`a la s\'erie
$\frac{1}{u^2} - \frac{a_1}{a_0 u} - \frac{a_0}{a(u) u^2}$.
L'identit\'e
$$
A_x\binom{0,1,\,\dots,\,k}{0,1,\,\dots,\,k} \; = \; A_{x=0}\binom{0,1,
\,\dots,\,k}{0,1,\,\dots,\,k}
$$
est d'ailleurs une illustration de la derni\`ere partie du
th\'eor\`eme \ref{degp0}.
\end{remark}

\noindent{\bf Preuve de l'assertion (ii) du th\'eor\`eme \ref{laymangen}.}
En appliquant l'assertion (i) du th\'eor\`eme \ref{laymangen}
\`a l'identit\'e $I^x[\frac{1}{t}\sum_{n=1}^\infty Q_n(0)t^n]=
\frac{1}{t}\sum_{n=1}^\infty Q_n(-x)t^n$, nous pouvons supposer $x=0$.
On a maintenant le d\'eveloppement
\begin{eqnarray*}
q_{s_1=0} (u)
&=&
s_0u + s_0^2s_2 u^3 + s_0^3s_3 u^4  + (s_0^4s_4 + 2s_0^3 s_2^2) u^5 + \cdots  \\
&=&s_0
\frac{u}{1-p(0)u-\displaystyle\frac{q(0)u^2}{1-p(1)u-\displaystyle\frac{q(1)u^2}{1-p(2)u-
\displaystyle\frac{q(2)u^2}{\ddots}}}}
\end{eqnarray*}
avec $p(0) = 0$, $q(0) = s_0s_2$, $p(1) =\frac{s_0s_3}{s_2}$,
$q(1) = s_0s_2+\frac{s_0^2s_4}{s_2}-\frac{s_0^2s_3^2}{s_2^2}$, \dots  \quad
La proposition \ref{indeps1} du chapitre \ref{Lukas} implique donc
\begin{eqnarray*}
q (u)
&=&
q_{s_1=0} \left(\frac{u}{1-s_1u}\right)\end{eqnarray*} 
et nous avons
\begin{eqnarray*}
q\left(\frac{u}{1-s_1u}\right)&=&s_0
\frac{u}{1-\bigl(s_1+p(0)\bigr)u-\displaystyle\frac{q(0)u^2}{1-\bigl(s_1+p(1)\bigr)u-
\displaystyle\frac{q(1)u^2}{1-\bigl(s_1+p(2)\bigr)u-\displaystyle\frac{q(2)u^2}{\ddots}}}}
\end{eqnarray*}
L'assertion (ii) d\'ecoule maintenant de la derni\`ere partie du th\'eor\`eme
\ref{degp0}.\hfill$\Box$

\begin{remark} L'identit\'e
$$\left(\frac{t}{1-tx}\right)^n=\sum_{k=n-1}^\infty \binom{k}{n-1}x^{k+1-n}\
t^{k+1}$$ montre qu'on a
$$\frac{1}{(1-xt)}\ a\left(\frac{t}{1-xt}\right)=\sum_{n,k} \binom{k}{n}
\ a_n\ x^{k-n}\ t^k$$
pour $a(t)=\sum_{n=0}^\infty a_n t^n$. La suite form\'ee des coefficients $b_k=
\sum_{n=0}^k \binom{k}{n}a_n$ est la {\it transform\'ee
binomiale} (de param\`etre $x$) de la suite $a=(a_0,\, a_1,\,\dots)$.
Il d\'ecoule de la
preuve ci-dessus que deux suites reli\'ees par une transformation
binomiale poss\`edent la m\^eme transform\'ee de Hankel.
\end{remark}

Les auteurs remercient Pierre de la Harpe et Fr\'ed\'eric Chapoton
pour des remarques et discussions int\'eressantes
ainsi que le Fonds National
Suisse de la Rechercher Scientifique pour un un soutien financier.


\vskip1cm

Roland Bacher

INSTITUT FOURIER

UMR 5582

BP 74

38402 St MARTIN D'H\`ERES Cedex (France)
 
e-mail: Roland.Bacher@ujf-grenoble.fr

\vskip1cm

Bodo Lass

INSTITUT DESARGUES

UMR 5028

21, Av. Claude Bernard

69622 VILLEURBANNE Cedex (France)
 
e-mail: lass@math.univ-lyon1.fr 


\begin{thebibliography}{10}

\bibitem{Baint} R.~Bacher, {\it Sur le groupe d'interpolation}, 
arXiv : math.CO/0609736.

\bibitem{Brom} Bromwich, An introduction to the theory of infinite
  series (Second edition revised with the assistance of
  T. M. Macrobert), St Martin's Press; Macmillan (1959).

\bibitem{C}
R.~Cori, \emph{Words and Trees},
Chapitre~11 dans le livre~{\cite{L}}.

\bibitem{D}
A.~Dvoretzky et Th.~Motzkin, \emph{A problem on arrangements},
Duke Math.~J., \textbf{14} (1947), 305-313.

\bibitem{F}
P.~Flajolet, \emph{Combinatorial aspects of continued fractions},
Discrete Math., \textbf{32} (1980), 125--161.

\bibitem{FZ} S. Fomin, A. Zelevinsky, {\it Total positivity: Tests and
parametrizations}, Math. Intell. {\bf 22}, No. 1 (2000), 23-33.

\bibitem{G}
J.~Gilewicz, \emph{Approximants de Pad\'e},
Lecture Notes in Mathematics, \textbf{667}, Springer (1978).

\bibitem{GJ} I.P. Goulden, J.M. Jackson, Combinatorial
enumeration, John Wiley \and Sons Ltd (1983).

\bibitem{Goursat} E. Goursat, Cours d'Analyse, Tome II, 7-i\`eme \'ed,
Gauthier-Villars (1949).

\bibitem{H} P. Henrici, Applied and computational complex analysis,
Volume I. Wiley Classics Library. New York etc:  John Wiley \and Sons Ltd
(1988).

\bibitem{H1}
P. Henrici, {\it Die Lagrange-B\"urmannsche Formel bei formalen Potenzreihen},
Jahresber. Deutsch. Math.-Verein.  {\bf 86} no. 4 (1984), 115-134.  

\bibitem{Kr} C. Krattenthaler, {\it Advanced determinant calculus},
S\'em. Loth. de Comb. 42 (1999), Article B42q.

\bibitem{Lagr} Lagrange: {\emph 
Nouvelle m\'ethode pour r\'esoudre des \'equations 
litt\'erales par le moyen des s\'eries}, M\'em. Acad. Roy. 
Belles-Lettres de Berlin {\bf XXIV} (1770) dans Oeuvres de Lagrange,
tome III, Gauthiers-Villars (1869), 5-73.

\bibitem{L}
M.~Lothaire, \emph{Combinatorics on Words},
Encyclopedia of Math. and its Applications, \textbf{17} (1983).

\bibitem{MOS} C. L. Mallows, A. M. Odlyzko, N. J. A. Sloane, {\it
Upper bounds for modular forms, lattices, and codes.}  J. Algebra
{\bf 36}  (1975), no. 1, 68--76.

\bibitem{Polya}
G. Polya, G. Szeg\"o, Problems and theorems in analysis I, 
Springer (1972).

\bibitem{R}
G.~N.~Raney, \emph{Functional composition patterns and power series reversion},
Trans.~Amer. Math.~Soc., \textbf{94} (1960), 441-451.

\bibitem{Stanl} R.~P.~Stanley, Enumerative Combinatorics, Volume 2, 
Cambridge University Press (1999).

\bibitem{V}
G.~Viennot, \emph{Une th\'eorie combinatoire des polyn\^omes orthogonaux g\'en\'eraux},
Notes de con\-f\'erences donn\'ees \`a l'Universit\'e du Qu\'ebec \`a Montr\'eal, 1983.

\bibitem{V1}
G.~Viennot, \emph{
A combinatorial theory for general orthogonal polynomials with
extensions and applications. Orthogonal polynomials and applications}, 
(Bar-le-Duc, 1984), Lecture Notes in Mathematics, {\bf 1171}, 
Springer (1985), 139-157. 

\bibitem{White} E.~T. Whittaker, G.~N. Watson, A course of modern
  analysis (4-th edition), Cambrige University Press (1978).

\end{thebibliography}
\end{document}